\documentclass[a4paper,10pt,fleqn]{article}
\usepackage{amsmath,amssymb,amsthm,graphicx,subfigure,float,caption,epstopdf,tabularx,color, bm,amsfonts,epic}
\usepackage[top=1.25in, bottom=1.0in, left=1.0in, right=1.0in]{geometry}
\usepackage{appendix}
\usepackage{amssymb}
\usepackage{multirow}
\usepackage{mathrsfs}
\usepackage[table]{xcolor}
\usepackage[numbers,sort&compress]{natbib}
\usepackage{comment,enumerate,multicol,xspace}
\usepackage{fancyhdr}
\usepackage{soul}
\usepackage{enumerate}
\usepackage{makecell}
\definecolor{Cobalt}{rgb}{0.25,0.41,0.88}
\usepackage[linkcolor=Cobalt,anchorcolor=Cobalt,citecolor=Cobalt,colorlinks]{hyperref}
\usepackage{bbm}
\usepackage{algorithm}
\usepackage{algorithmic}
\usepackage{booktabs}
\usepackage{lineno}
\usepackage{dsfont}
\usepackage{pifont}

\allowdisplaybreaks
\newtheorem{thm}{Theorem}[section]
\newtheorem{lem}{Lemma}[section]

\numberwithin{equation}{section}

\usepackage{indentfirst}
\graphicspath{{Fig}}
\begin{document}
\title{Finite element method with Gr\"unwald-Letnikov type approximation in time for a constant time delay subdiffusion equation}

\author{Weiping Bu$^{1,2}$\footnote{Correspondence author. Email: weipingbu@xtu.edu.cn},\quad Xueqin Zhang$^1$,\quad Weizhi Liao$^1$,\quad Yue Zhao$^3$
\\{\small $^1$ School of Mathematics and Computational Science, Xiangtan University, }
\\{\small Hunan 411105, China}
\\{\small $^2$ Hunan Key Laboratory for Computation and Simulation in Science and}
\\{\small Engineering, Hunan 411105, China}
\\{\small$^3$ Space Engineering University, Beijing 101400, China}
}

\date{}
\maketitle
\begin{abstract}
In this work, a subdiffusion equation with constant time
delay $\tau$ is considered.
First, the regularity of the solution to the
considered problem is investigated, finding that its first-order time derivative exhibits
singularity at $t=0^+$ and its second-order time derivative shows singularity at both
$t=0^+$ and $\tau^+$, while the solution can be decomposed
into its singular and regular components. Then, we derive a fully discrete
finite element scheme to solve the considered problem based on the standard Galerkin finite element
method in space and the Gr\"unwald-Letnikov type approximation in
time. The analysis shows that the developed numerical scheme is stable.
In order to discuss the error estimate, a new discrete Gronwall inequality is established. Under the above decomposition
of the solution, we obtain a local error estimate in time for the developed numerical scheme. Finally,
some numerical tests are provided to support our theoretical analysis.
\\[2ex]
\textbf{AMS subject classification:} 65M06, 65M12, 65M60.\\[2ex]
\textbf{Keywords:} Subdiffusion equation with constant time delay; finite element method;  Gr\"unwald-Letnikov type approximation; discrete Gronwall inequality; error estimate
\end{abstract}

%%%%%%%%%%%%%%%%%%%%%%%%%%%%%%%%%%%%%%%%%% end abstract %%%%%%%%%%%%%%%%%%%%%%%%%%%%%%%%%%%%%%%%%%%%%%%%%%%%%%%%%%%

 \vskip 5mm
 %%%%%%%%%%%%%%%%%%%%%%%%%%%%%%%%%%%%%%%%%% begin section 1 introduction %%%%%%%%%%%%%%%%%%%%%%%%%%%%%%%%%%%%%%%%%%%%
\section{Introduction}
\label{Sec1}
It is well known that differential equation is a powerful tool for understanding and describing some
complex phenomena in the natural world.
In view of the fact that there exist numerous real-world systems which depend on their historical states, this makes the delay differential equation show wide application potential in various fields such as physics, biology and medicine \cite{Erneux2009applied,Smith2011An}.
The emergence of this trend has aroused great interest of mathematicians in delay differential equations \cite{Driver1977ordinary,Gopalsamy1992Stability,Bellen2003numerical,Balachandran2009delay}.
Compared with integer-order model, fractional differential equation often exhibits high
reliability in capturing memory effect, hereditary characteristic and non-locality
of complex systems. Therefore, the fractional models often are used to describe these
practical problems \cite{Carpinteri2014Fractals,Hilfer2000Applications,Kilbas2006Theory,Sun2018A}.
In this paper, we consider the following time-fractional diffusion equation with constant time delay
\begin{equation}\label{eqs1_1}
\partial_{t} u(x,t)=\partial_{t}^{1-\alpha}(p\Delta u(x,t)+au(x,t))+bu(x,t-\tau)+f(x,t),\quad
(x,t)\in \Omega \times (0,K\tau],
\end{equation}
\begin{equation}\label{eqs1_2}
u(x,t)=\varphi(x,t),\quad (x,t)\in \bar{\Omega} \times \left[-\tau,0\right],
\end{equation}
\begin{equation}\label{eqs1_3}
 u(x,t)|_{\partial \Omega}=0, \quad t \in [-\tau,K\tau],
\end{equation}
where $\Omega = (0, L)$, $L$,$p,a,b$ are some constants satisfying $L\geq0,p\geq0,a\leq0$,$b\neq 0$,
$f(x,t)$ are continuous functions, $\tau>0$ is the time delay parameter, $K$ is a given
positive integer, $\varphi(x,t)$ is a continuous function in $\bar{\Omega} \times \left[-\tau,0\right]$, and $\partial_t^lu$ denote as the $l$th-order derivative with respect to variable $t$ and $\partial_t(\cdot):=\partial^{1}_t(\cdot)$ . The Riemann-Liouville fractional derivative $\partial_{t}^{1-\alpha}u(x,t)$ is defined by \cite{Bu2024Finite}
$$\partial_{t}^{1-\alpha}u(x,t)=\frac{\partial}{\partial t}\int_0^t\omega_\alpha(t-s)u(x,s)ds,\quad 0<\alpha<1,$$ where $\omega_\beta(t):=\frac{t^{\beta-1}}{\Gamma(\beta)}$ for $\beta, t>0$ is the convolution kernel function.

At present, there have been some advancements in the theoretical studies of fractional
diffusion equations with delay. In \cite{Prakash2020Exact}, Prakash et al.
proposed the invariant subspace approach to find the exact solution of time-fractional
reaction-diffusion equation with time delay. In \cite{Zhu2016Local}, Zhu et al.
investigated the local and nonlocal existence of mild solution for a
nonlinear time delay fractional reaction-diffusion equation.
By using the semigroup theory of operators and the monotone iterative technique, Ref.
\cite{Li2021Monotone} obtained the existence and uniqueness of mild solutions for a
time-space diffusion equation involving delay. Yao and Yang \cite{Yao2023Stability} investigated the
asymptotic stability and long-time decay rates for a fractional diffusion-wave equation
with time delay. Meanwhile, some researchers proposed a lot of effective numerical methods
to solve the time fractional diffusion equation with delay. For example, in \cite{Li2018convergence,zhao2018fast}, researchers
discussed finite difference methods for constant time delay subdiffusion
equation. In \cite{Li2018ane,Zhang2023numerical,Alikhanov2024second},
some finite difference schemes were devised to solve the time fractional
diffusion-wave equations with delay. We also noted that some studiers use
finite difference method to solve distributed order fractional diffusion equation
\cite{Pimenov2017on}. In \cite{peng2023convergence,Peng2024uncond}, the authors
considered the finite element methods for the time delay subdiffusion equation and
investigated their convergence and superconvergence. Based on L1 formula in time and
Galerkin spectral method in space, Zaky et al. \cite{zaky2023L1} developed an
effective numerical method
for variable-order time fractional reaction-diffusion equation with delay.
In order to analyze the stability and convergence of finite difference scheme of time
factional multi-delayed diffusion equation,
Hendy and Mac\'ias-D\'iaz \cite{Hendy2019novel} proposed a
novel discrete Gronwall inequality.

Recently, Tan et al. \cite{Tan2022L1} discuss a constant time delay fractional diffusion
equation from both theoretical and numerical perspectives. They not only derive its exact
solution, but also reveal the multiple singularity phenomena of the solution in time.
Furthermore, Cen and Vong \cite{cen2023the} obtain a sharp multi-singularity
result of a class of delay fractional ordinary differential equation, and two
corrected L-type schemes are constructed in \cite{cen2023corrected}
to overcome this multi-singularity.
Subsequently, by using variable-step L1 method,
Refs. \cite{Bu2024Finite,Ou2024Variable} obtain effective numerical schemes with
satisfactory temporal convergence accuracy when the exact solution has multiple
singularities in time.
%It is worth pointing out that the current work on
%numerically solving the time-fractional delay diffusion equation based on multi-point singularity is still
%limited, which inspires us to study this problem further.
It is worth pointing out that the regularity results of Refs.
\cite{Bu2024Finite,cen2023the,Ou2024Variable} show that
the first derivative of the solution with respective to variable $t$ satisfies
\begin{equation}\label{eqs1_1207a}
|\partial_t u|\leq C\left(1+\left(t-(k-1)\tau\right)^{k\alpha-1}\right), \ t\in ((k-1)\tau,k\tau],
\end{equation}
where $k=1,2,\cdots.$ However, compared with previous works, we discover that there are
some different regularity results for the present problem (\ref{eqs1_1})--(\ref{eqs1_3}).
In this paper,
the first derivative of the solution with respective to variable $t$
has
\begin{equation}\label{eqs1_1207b}
|\partial_t u|\leq C(1+t^{\alpha-1}), \ t\in ((k-1)\tau,k\tau],
\end{equation}
and the second derivative of the solution with respective to variable $t$ satisfies
\begin{equation}\label{eqs1_1207c}
|\partial_t^2 u|\leq C(1+t^{\alpha-2}), \ t\in(0,\tau]\ \text{and}\
|\partial_t^2 u|\leq C(1+(t-\tau)^{\alpha-1}), \ t\in((k-1)\tau,k\tau],
\end{equation}
where $k=2,3,\cdots.$ The regularity results (\ref{eqs1_1207b}) and
(\ref{eqs1_1207c}) indicate that $\partial_t u$ will blow up only at $t=0$,
and $\partial_t^2 u$ will blow up only at $t=0$ and $t=\tau$. It implies that
the problem (\ref{eqs1_1})--(\ref{eqs1_3}) has better regularity than the problem
considered in \cite{Bu2024Finite,cen2023the,Ou2024Variable}.
What is more, the solution can be decomposed
into its singular and regular components.
Therefore, it motivates us
to develop new numerical method to handle the different
smoothness.
In view of the good regularity of exact solution and in order
to avoid the round-off error which may be caused by dense graded mesh, we use the
Gr\"unwald-Letnikov approximation to discretize the time fractional operators, and propose
a fully discrete finite element scheme for (\ref{eqs1_1})--(\ref{eqs1_3}). Then we
discuss the stability of the developed numerical scheme, and investigate the convergence
based on the regularity results (\ref{eqs1_1207b}) and
(\ref{eqs1_1207c}).

The structure of this work is as follows. In Section \ref{Sec2}, the exact solution of
(\ref{eqs1_1})--(\ref{eqs1_3}) is obtained and its regularity is discussed. In Section \ref{Sec3}, we establish a fully discrete finite element scheme for the problem
(\ref{eqs1_1})--(\ref{eqs1_3}), and investigate the stability. Based on the regularity results, the convergence is discussed in Section \ref{Sec4}.
Finally, numerical tests are presented to verify the theoretical results in Section \ref{Sec5} and the conclusion is given in Section \ref{Sec6}.

\textbf{Notation.} In this paper, $C$ is a general constant, which can be different in
different situations. Denote $(\cdot,\cdot)$ and $\|\cdot\|_0$ as the inner product and
norm of $L^2(0,L)$. Define the solution $u(x,t),t\in((k-1)\tau,k\tau]$ of (\ref{eqs1_1})--(\ref{eqs1_3}) by $u_{k\tau}(x,t)$, and $u_{\tau}(x,t):=u_{1\tau}(x,t)$.

\section{The regularity of the solution}\label{Sec2}
In this section, we discuss the regularity of the solution to (\ref{eqs1_1})--(\ref{eqs1_3})
by the variable separation method and the Laplace transform method.
First, when $(x,t)\in[0,L] \times (0,\tau]$, the problem (\ref{eqs1_1})--(\ref{eqs1_3}) can be written as
\begin{equation}\label{eqs2_1}
	\left\{ \begin{aligned}
		\partial_{t}u(x,t)&=\partial_{t}^{1-\alpha}(-G(u))+F_{\tau}(x,t),\ &0<x<L,\ &0<t\leq \tau,\\
		u(x,t)&=\varphi(x,t),\ &0\leq x\leq L,&-\tau\leq t\leq 0,\\
		u(0,t)&=0,u(L,t)=0,\ & -\tau\leq t\leq \tau ,
	\end{aligned} \right.
\end{equation}
where $G(u):=-p\Delta u(x,t)-au(x,t),\ F_{\tau}(x,t):=bu(x,t-\tau)+f(x,t)$.
Similar to \cite{Bu2024Finite}, it follows from the variable separation method that
\begin{equation}\label{eqs2_2}
	u_{\tau}(x,t)=\sum_{i=1}^{\infty}T_{i}(t)X_{i}(x),\quad T_{i}(t)=(u_{\tau}(x,t),X_{i}),
\end{equation}
where
$X_{i}$ is the eigenfunction corresponding to the eigenvalue $\lambda_i >0$ satifying $\Vert X_i \Vert_0=1$ and
\begin{equation}G(X_{i})=\lambda_{i}X_{i},\quad X_i(0) = X_i(L) = 0,
\end{equation}
and
\begin{equation}\label{eqs2_4}
\varphi(x,t-\tau)=\sum_{i=1}^{\infty} \varphi_i(t-\tau) X_i(x), \quad F_\tau(x, t)=\sum_{i=1}^{\infty} F_{\tau, i}(t) X_i(x),
\end{equation}
where
$\varphi_i(x,t-\tau)=(\varphi(x,t-\tau),X_i),\ F_{\tau, i}=(F_{\tau}(x,t),X_i), t\in(0,\tau]$.
Combining (\ref{eqs2_1}), (\ref{eqs2_2}) and (\ref{eqs2_4}), then we have
\begin{align*}
    \partial_{t}T_{i}(t)=\partial_{t}^{1-\alpha}(-\lambda_{i}T_{i}(t))+F_{\tau,i}(t).
\end{align*}
Taking the Laplace transform of the above equation yields
\begin{align*}	s\tilde{T_{i}}(s)-\varphi_{i}(0)
=-s^{1-\alpha}\lambda_{i}\tilde{T_{i}}(s)+\tilde{F}_{\tau,i}(s),	
\end{align*}
i.e.,
\begin{align*} \tilde{T_{i}}(s)=\frac{\varphi_{i}(0)}{s+\lambda_{i}s^{1-\alpha}}+\frac{\tilde{F}_{\tau,i}(s)}{s+\lambda_{i}s^{1-\alpha}}.
\end{align*}
Furthermore, from the inverse Laplace transform, one has
\begin{equation}\label{eqs2_1208c}
T_{i}(t)=\varphi_{i}(0)E_{\alpha}(-\lambda_{i}t^{\alpha})+\int_{0}^{t}E_{\alpha}(-\lambda_{i}s^{\alpha})F_{\tau,i}(t-s)ds,
\end{equation}
where $E_{\mu,\nu}(z):=\sum_{j=0}^{\infty}
\frac{z^j}{\Gamma(j\mu+\nu)}$ is the Mittag-Leffler function.
Therefore, we can obtain the solution of (\ref{eqs2_1}) by substituting the above equation into (\ref{eqs2_2}) as follows
\begin{equation}\label{eqs2_5}
u_{\tau}(x,t)=\sum_{i=1}^{\infty}\left(\varphi_{i}(0)E_{\alpha}(-\lambda_{i}t^{\alpha})+\int_{0}^{t}E_{\alpha}(-\lambda_{i}s^{\alpha})
F_{\tau,i}(t-s)ds\right)X_{i}(x).
\end{equation}
In fact, (\ref{eqs2_5}) implies that the problem (\ref{eqs1_1})--(\ref{eqs1_3})
exists a unique solution in the domain $[0,L]\times (0,\tau]$. Moreover,
by performing similar manipulation to $[0,L]\times ((k-1)\tau,k\tau], k=2,3,\dots,K$, it is clear that the problem (\ref{eqs1_1})--(\ref{eqs1_3})
has a unique solution in $[0,L]\times ((k-1)\tau,k\tau]$, i.e.,
\begin{equation}\label{eqs2_6}
u_{k\tau}(x,t)=\sum_{i=1}^{\infty}\left(\varphi_{i}(0)E_{\alpha}(-\lambda_{i}t^{\alpha})+\int_{0}^{t}E_{\alpha}(-\lambda_{i}s^{\alpha})
F_{k\tau,i}(t-s)ds\right)X_{i}(x),
\end{equation}
where
$\varphi_{i}(0)=\left(\varphi(x,0),X_{i}\right),\ F_{k\tau,i}(t)=(F_{k\tau}(x,t),X_{i}(x))$ and
\begin{equation}\label{eqs2_7}
    \begin{aligned}
     F_{k\tau}(x,t)=\left\{ \begin{aligned}
		&b\varphi(x,t-\tau)+f(x,t),&0\leq x\leq L,\ &0\leq t\leq \tau,\\
		&bu_{\tau}(x,t-\tau)+f(x,t),&0\leq x\leq L,\ &\tau<t\leq 2\tau,\\
        &\vdots\\
		&bu_{(k-1)\tau}(x,t-\tau)+f(x,t),\ &0\leq x\leq L,&(k-1)\tau<t\leq k\tau.
	\end{aligned} \right.
    \end{aligned}
\end{equation}

In order to investigate the regularity of the solution, we give two
auxiliary lemmas firstly.
\begin{lem}\label{jullemma2.1}(\cite{Sakamoto2011Initial})
If $\alpha>0,\lambda>0,t>0$ and positive integer  $m\in \mathbb{N}$, then
$$
\frac{d^m}{dt^m}E_{\alpha,1}(-\lambda t^\alpha)=-\lambda t^{\alpha-m}E_{\alpha,\alpha-m+1}(-\lambda t^\alpha).$$
Moreover, if $0<\alpha<1$ and $\eta \geq 0$, then
$E_{\alpha,\alpha}(-\eta)\geq 0.$
\end{lem}

From Lemma \ref{jullemma2.1}, it is easy to know that
\begin{equation}\label{eqs2_1208a}
\int_0^t|\lambda_is^{\alpha-1}E_{\alpha,\alpha}(-\lambda_is^\alpha)|ds =1-E_{\alpha,1}(-\lambda_it^\alpha)\leq C.
\end{equation}
Let $\gamma>0.$ Define the space $D\left(G^\gamma\right)$ and its norm as follows \cite{Bu2024Finite}
$$
D\left(G^\gamma\right)=\left\{g \in L^2(0, L): \sum_{i=1}^{\infty} \lambda_i^{2 \gamma}\left|\left(X_i, g\right)\right|^2\leq \infty\right\},
$$
and
$$
\|g\|_{G^\gamma}=\left(\sum_{i=1}^{\infty} \lambda_i^{2 \gamma}\left|\left(X_i, g\right)\right|^2\right)^{1/2}.
$$
\begin{lem}\label{jullemma2.2}
Suppose that\\
$\mathrm{(I)}\;\left\|\varphi(x,t-\tau)\right\|_{G^{3/2}}\leq C$ and $\left\|\partial_{t} \varphi(x,t-\tau)\right\|_{G^{1/2}}\leq C(1+t^{\alpha-1})$, $t\in[0,\tau]$;\\
$\mathrm{(II)}\; \left\| f(x,t)\right\|_{G^{3/2}}\leq C,\ t\in[0,K\tau]$, $\left\|\partial_{t} f(x,t)\right\|_{G^{1/2}}\leq C$, $t\in(0,\tau]$ and $\left\|\partial_{t} f(x,t)\right\|_{G^{1/2}}\leq C(1+(t-\tau)^{\alpha-1})$, $t\in(\tau,K\tau]$\\
Then for $k\geq 1,$ $\|u_{k\tau }(x,t)\|_{G^{3/2}}, \|F_{k\tau}(x,t)\|_{G^{3/2}}\in L^{\infty}\left(((k-1)\tau ,k\tau]\right), \|
\partial_{t}u_{k\tau}(x,t)\|^2_{G^{1/2}}
\leq C\left(1+t^{2(\alpha-1)}\right)$, $\|\partial_tF_{\tau}(x,t)\|_{G^{1/2}}^2\leq C$, and for $k\geq2,$ $t\in (j\tau,(j+1)\tau]$ with $0\leq j\leq k-1$, $\|\partial_tF_{k\tau}(x,t)\|_{G^{1/2}}^2\leq C\left(1+(t-\tau)^{2(\alpha-1)}\right)$.
%\begin{align*}
%   \begin{aligned}
%&\|\partial_tF_{k\tau}(t)\|_{G^{1/2}}^2\leq C\ \text{for}\ k=1,\\
%&\|\partial_tF_{k\tau}(t)\|_{G^{1/2}}^2\leq C\left[1+(t-\tau)^{2(\alpha-1)}\right]\ \text{for}\ k=2,3,\cdots,K.
%\end{aligned}
%\end{align*}
\end{lem}
\noindent\textbf{Proof.}
    When $k=1$, it follows from (\ref{eqs2_7}), the conditions (I) and (II) that
\begin{equation}\label{eqs2_1208b}
\begin{aligned}
\left\|F_\tau(x,t)\right\|^2_{G^{3/2}}=&\left\|b\varphi (x,t-\tau)+f (x,t)\right\|^2_{G^{3/2}}\\
\leq& C\left[\left\|\varphi (x,t-\tau)\right\|^2_{G^{3/2}}+\left\|f (x,t)\right\|^2_{G^{3/2}}\right].
\end{aligned}
\end{equation}
Therefore, (\ref{eqs2_1208b}) implies that $\|F_\tau(x,t)\|_{G^{3/2}}\in L^{\infty}([0,\tau])$.
By (\ref{eqs2_1208c}) and the boundedness of Mittag-Leffler function, we know that
\begin{equation}\label{eqs2_1208d}
\begin{aligned}
|T_i(t)|^2\leq&\left(|\varphi_i(0)|+\int_0^t|F_{\tau,i}(t-s)|ds\right)^2\\
\leq&C\left[|\varphi_i(0)|^2+\int_0^t|F_{\tau,i}(t-s)|^2ds \right].
\end{aligned}
\end{equation}
The combination of (\ref{eqs2_1208c}), (\ref{eqs2_5}) and (\ref{eqs2_1208d})
yields
\begin{equation}\label{eqs2_1208e}
\begin{aligned}
\|u_\tau(x,t)\|_{G^{3/2}}^2=&\sum_{i=1}^{\infty}\lambda_i^3|T_i(t)|^2\\
\leq&C\left[\|\varphi(x,0)\|_{G^{3/2}}^2+\int_0^{\tau}\left\|F_\tau(x,t-s)\right\|_{G^{3/2}}^2ds \right]\leq C.
\end{aligned}
\end{equation}
Hence, (\ref{eqs2_1208e}) indicates that $\|u_{\tau}(x,t)\|_{G^{3/2}}\in L^{\infty}((0,\tau])$.
Now, we suppose that $\|u_{k\tau}(x,t)\|_{G^{3/2}}, \|F_{k\tau}(x,t)\|_{G^{3/2}}\in L^{\infty}(((k-1)\tau ,k\tau])$ are true for $k\leq n$. When $k=n+1$, (\ref{eqs2_7}) gives
\begin{equation}\label{eqs2_1208f}
\begin{aligned}
\|F_{(n+1)\tau}(x,t)\|^2_{G^{3/2}}=&\|bu_{n\tau} (x,t-\tau)+f (x,t)\|^2_{G^{3/2}}\\
\leq& C\left[\|u_{n\tau} (x,t-\tau)\|^2_{G^{3/2}}+\|f (x,t)\|^2_{G^{3/2}}\right]\leq C.
\end{aligned}
\end{equation}
Besides, from (\ref{eqs2_6}), one has
$$
\begin{aligned}
|u_{(n+1)\tau,i}(t)|^2&\leq C\left[|\varphi_i(0)|^2+\int_0^{(n+1)\tau}|F_{(n+1)\tau,i}(t-s)|^2ds\right].
\end{aligned}
$$
The above inequality leads to
\begin{equation}\label{eqs2_1208g}
\begin{aligned}
     \|u_{(n+1)\tau}(x,t)\|_{G^{3/2}}^2\leq C\left[\|\varphi(x,0)\|_{G^{3/2}}^2+\int_0^{(n+1)\tau}\|F_{(n+1)\tau}(x,s)\|_{G^{3/2}}^2ds\right]\leq C.
\end{aligned}
\end{equation}
According to (\ref{eqs2_1208f}), (\ref{eqs2_1208g}) and the mathematical induction, we can conclude that $\|u_{k\tau}(x,t)\|_{G^{3/2}}, \|F_{k\tau}(x,t)\|_{G^{3/2}}\in L^{\infty}(((k-1)\tau ,k\tau])$ hold for $k=1,2,\cdots,K$.

Next, we discuss $\|\partial_{t}u_{k\tau}(t)\|_{G^{1/2}}$ and $\|\partial_tF_{k\tau}(t)\|_{G^{1/2}}.$ First, it is easy to know that the above results imply
\begin{equation}\label{eqs2_1208h}
 \begin{aligned}
\sum_{i=1}^{\infty}|F_{k\tau,i}(t)|\leq \left(\sum_{i=1}^{\infty}\frac{1}{\lambda_i}\right)^{\frac{1}{2}}\left(\sum_{i=1}^{\infty}\lambda_i|F_{k\tau,i}(t)|^2\right)^{\frac{1}{2}}  \leq C,
 \end{aligned}
\end{equation}
and
\begin{equation}\label{eqs2_1209a}
    \begin{aligned}
\sum_{i=1}^{\infty}\lambda_i|F_{k\tau,i}(t)|\leq \left(\sum_{i=1}^{\infty}\frac{1}{\lambda_i}\right)^{\frac{1}{2}}\left(\sum_{i=1}^{\infty}\lambda^3_i|F_{k\tau,i}(t)|^2\right)^{\frac{1}{2}}  \leq C,
    \end{aligned}
\end{equation}
where the well known fact $\lambda_i\approx i^2$ is used.
Taking the first derivative of $u_{k\tau} (x,t),\ k=1,2,\cdots,K$ with respect to $t$ and applying Lemma \ref{jullemma2.1}, it yields
\begin{equation}\label{eqs2_12}
  \begin{aligned} \partial_t u_{k\tau}(x, t)
 =&\sum_{i=1}^{\infty}\Big\{-\lambda_i t^{\alpha-1}\varphi_{i}(0)E_{\alpha,\alpha}(-\lambda_{i}t^{\alpha})+E_{\alpha}(-\lambda_{i}t^{\alpha})F_{k\tau,i}(0)\\ &+\int_{0}^{t}E_{\alpha}(-\lambda_{i}s^{\alpha})
\frac{d}{dt}F_{k\tau,i}(t-s)ds\Big\} X_{i}(x)\\=&\sum_{i=1}^{\infty}\Big\{-\lambda_i t^{\alpha-1}\varphi_{i}(0)E_{\alpha,\alpha}(-\lambda_{i}t^{\alpha})+E_{\alpha}(0)F_{k\tau,i}(t)\\&+\int_{0}^{t}(-\lambda_i  s^{\alpha-1})E_{\alpha,\alpha}(-\lambda_{i}s^{\alpha})
F_{k\tau,i}(t-s)ds\Big\} X_{i}(x),
\end{aligned}
\end{equation}
where the integration by parts is used. Therefore, one has
\begin{equation}\label{eqs2_17}
\begin{aligned}
|\partial_{t}u_{k\tau,i}(t)|^2&\leq C\left[t^{2(\alpha-1)}\lambda_i^2| \varphi_{i}(0)|^2+|F_{k\tau,i}(t)|^2
+\left(\int_{0}^{t}\lambda_i (t-s)^{\alpha-1}E_{\alpha,\alpha}(-\lambda_{i}(t-s)^{\alpha})
F_{k\tau,i}(s)ds\right)^2\right].
\end{aligned}
\end{equation}
Applying the Young's inequality for convolution, the Cauchy-Schwarz inequality and (\ref{eqs2_1208a}), it holds
$$
\begin{aligned}
|\partial_{t}u_{k\tau,i}(t)|^2
\leq& C\left[t^{2(\alpha-1)}\lambda_i^2 |\varphi_{i}(0)|^2+|F_{k\tau,i}(t)|^2+\left(\int_{0}^{t}| \lambda_is^{\alpha-1}E_{\alpha,\alpha}(-\lambda_{i}s^{\alpha})|ds\right)^2
\left(\int_{0}^{t}|
F_{k\tau,i}(s)|^2ds\right)\right]\\
\leq& C\left[t^{2(\alpha-1)}\lambda_i^2| \varphi_{i}(0)|^2+|F_{k\tau,i}(t)|^2
+
\int_{0}^{t}|
F_{k\tau,i}(s)|^2ds\right].
\end{aligned}
$$
The above inequality demonstrates
\begin{equation}
 \begin{aligned}
\|
\partial_{t}u_{k\tau}(x,t)\|^2_{G^{1/2}}
&\leq C\left[t^{2(\alpha-1)}\sum_{i=1}^{\infty}\lambda^3_i|\varphi_{i}(0)|^2+\sum_{i=1}^{\infty}\lambda_i|F_{k\tau,i}(t)|^2
+
\int_{0}^{t}\sum_{i=1}^{\infty}\lambda_i|
F_{k\tau,i}(s)|^2ds\right]\\
&\leq C\left(1+t^{2(\alpha-1)}\right).
\end{aligned}
\end{equation}
Since $\partial_tF_{k\tau}(x,t)=b\partial_t\varphi(x,t-\tau)+\partial_tf(x,t)$ for $k=1,$ and $\partial_tF_{k\tau}(x,t)=b\partial_tu_{j\tau}(x,t-\tau)+\partial_tf(x,t)$, $j\tau<t\leq (j+1)\tau$, $1\leq j\leq k-1$ for $k=2,3,\cdots,K,$
it is clear that the condition (I) and (II) give
\begin{equation}\label{eqs2_1209b}
 \begin{aligned}
\|\partial_tF_{k\tau}(x,t)\|_{G^{1/2}}^2\leq&C\left[\|\partial_t \varphi(x,t-\tau) \|_{G^{1/2}}^2 + \|\partial_t f(x,t) \|_{G^{1/2}}^2\right] \leq C\ \text{for}\ k=1,
\end{aligned}
\end{equation}
and for $k\geq 2$,
$1\leq j\leq k-1$, $t\in (j\tau,(j+1)\tau]$, one has
\begin{equation}\label{eqs2_1209c}
\begin{aligned}
\|\partial_tF_{k\tau}(x,t)\|_{G^{1/2}}^2&\leq C\left[\| \partial_tu_{j\tau}(x,t-\tau) \|_{G^{1/2}}^2 + \| \partial_tf(x,t)\|_{G^{1/2}}^2\right]\\
&\leq C\left(1+(t-\tau)^{2(\alpha-1)}\right).
\end{aligned}
\end{equation}
This complete the proof of Lemma \ref{jullemma2.2}.
$\quad\Box$

Now we state the regularity result.
\begin{thm}\label{thm2.1}
Suppose that the assumptions (I) and (II) of Lemma \ref{jullemma2.2} hold.
%$\mathrm{(I)}\;\varphi(x,t-\tau)\in D\left(G^{3 / 2}\right)$ and $\partial_{t} \varphi(x,t-\tau)\in D\left(G^{1 / 2}\right)$, $t\in[0,\tau]$;\\
%$\mathrm{(II)}\; f(x,t)\in D\left(G^{3 / 2}\right),\ \partial_{t} f(x,t)\in D\left(G^{1 / 2}\right)$ and $\int_0^t\|f(x,s)\|_{G^{3/2}}^2ds<\infty$, $t\in(0,K\tau]$.
Then the solution of (\ref{eqs1_1})--(\ref{eqs1_3}) which is expressed by (\ref{eqs2_5}) and (\ref{eqs2_6}) satisfies
$$
\begin{aligned}	&|\partial_{t}^{l}u_{k\tau}(x,t)|\leq C(1+t^{\alpha-l}),\ k=1,2,\cdots,K,\  l=0,1,\\
 &|\partial_{t}^{2}u_{\tau}(x,t)|\leq C(1+t^{\alpha-2}),\\
&|\partial_{t}^{2}u_{k\tau}(x,t)|\leq C(1+(t-\tau)^{\alpha-1}), \  k=2,3,\dots ,K,\\
 &|\partial_{x}^{n}u_{k\tau}(x,t)|\leq C,\ k=1,2,\cdots,K,\ n=0,1,2.
\end{aligned}
$$
%where $x\in [0,L].$
Furthermore, under the conditions (I) and (II) of Lemma \ref{jullemma2.2} and   $\left|G^{m+1}(\varphi(x,0))
 \right|\leq C$, the solution of (\ref{eqs1_1})--(\ref{eqs1_3}) can be decomposed as
\begin{equation}\label{eqs4_1}
    \begin{aligned}
    u(x,t)=\sum_{j=0}^m \gamma_jt^{j\alpha}+\tilde{\gamma}t+Y(x,t),\ (x,t)\in [0,L]\times[0,K\tau],
    \end{aligned}
\end{equation}
where $m$ is the smallest positive integer to ensure that $(m+1)\alpha> 1$, $\gamma_j$ and $\tilde{\gamma}$ are bounded functions with respect to the variable $x$, and $Y(\cdot,t)$ satisfies
\begin{equation}\label{eqs4_2}
    \begin{aligned}
Y(\cdot,0)=\partial_t Y(\cdot,t)|_{t=0}=0,\ Y(\cdot,t)\in C^1[0,K\tau],\ \left|\partial_tY(\cdot,t)\right|\leq C(1+t^{(m+1)\alpha-1}), \ \partial^2_tY(\cdot,t)\in L^1[0,K\tau].
    \end{aligned}
\end{equation}
\end{thm}
\noindent\textbf{Proof.}
First,  according to (\ref{eqs2_5}), (\ref{eqs2_6}) and (\ref{eqs2_1208h}), it is obvious that
\begin{equation}
\begin{aligned}
|u_{k\tau}(x,t)|&=\left|\sum_{i=1}^{\infty}\left[\varphi_{i}(0)E_{\alpha}(-\lambda_{i}t^{\alpha})+\int_{0}^{t}E_{\alpha}(-\lambda_{i}s^{\alpha})
F_{k\tau,i}(t-s)ds\right ]X_{i}(x)\right|
\\&\leq C\left[\sum_{i=1}^{\infty}|\varphi_{i}(0)|+\int_{0}^{t}\sum_{i=1}^{\infty}\left|
F_{k\tau,i}(t-s)\right|ds\right]
\\&\leq C(1+t^\alpha).
\end{aligned}
\end{equation}
Moreover, it follows from the condition (I) and (\ref{eqs2_1208h})--(\ref{eqs2_12}) that
\begin{equation}\label{eqs0311c}
  \begin{aligned}
    |\partial_t u_{k\tau}(x, t)|\leq & C\left[t^{\alpha-1}\sum_{i=1}^{\infty}\lambda_i |\varphi_{i}(0)|+\sum_{i=1}^{\infty}|F_{k\tau,i}(t)|+\int_0^t(t-s)^{\alpha-1}\sum_{i=1}^{\infty}\lambda_i|F_{k\tau,i}(s)|ds\right] \\\leq& C(1+t^{\alpha-1}).
\end{aligned}
\end{equation}
Now we discuss $|\partial_{t}^2u_{k\tau}(x,t)|,\ k=1,2,\cdots,K.$
Taking the first time derivative on both sides of
(\ref{eqs2_12}) yields
\begin{equation}\label{eqs0315a}
\begin{aligned}
\partial_t^2 u_{k\tau}(x,t) =&\sum_{i=1}^{\infty}\Big\{-\lambda_i t^{\alpha-2}E_{\alpha,\alpha-1}(-\lambda_it^\alpha)\varphi_{i}(0)+E_{\alpha}(0)\partial_tF_{k\tau,i}(t)\\
 &+(-\lambda_{i}t^{\alpha-1}E_{\alpha,\alpha}(-\lambda_{i}t^{\alpha}))F_{k\tau,i}(0)+\int_{0}^{t}-\lambda_i s^{\alpha-1}E_{\alpha,\alpha}(-\lambda_{i}s^{\alpha})
\partial_tF_{k\tau,i}(t-s)ds \Big\}X_i(x).
\end{aligned}
\end{equation}
Furthermore, we have
\begin{equation}\label{eqs2_24}
    \begin{aligned}
    |\partial_t^2 u_{k\tau}(x,t)|\leq& C\left[t^{\alpha-2}\sum_{i=1}^{\infty}\lambda_i|\varphi_{i}(0)|+\sum_{i=1}^{\infty}\left|\partial_tF_{k\tau,i}(t)\right|\right.\\
 &\left.+t^{\alpha-1}\sum_{i=1}^{\infty}\lambda_i|F_{k\tau,i}(0)|+\sum_{i=1}^{\infty}
\left|\int_{0}^{t}-\lambda_i s^{\alpha-1}E_{\alpha,\alpha}(-\lambda_{i}s^{\alpha})
\partial_tF_{k\tau,i}(t-s)ds\right|\right]\\
\leq& C\left[t^{\alpha-2}\|\varphi(x,0)\|_{G^{3/2}}+\|\partial_tF_{k\tau}(x,t)\|_{G^{1/2}}+t^{\alpha-1}\|F_{k\tau}(x,0)\|_{G^{3/2}}+\int_{0}^{t}\|\partial_sF_{k\tau}(x,s)\|_{G^{1/2}}ds\right].\\
\end{aligned}
\end{equation}
According to Lemma \ref{jullemma2.2},
when $k=1$ it yields
\begin{equation}\label{eqs0311a}
  \begin{aligned}
    |\partial_t^2 u_{\tau}(x,t)|
\leq C\left(1+t^{\alpha-2}\right),
\end{aligned}
\end{equation}
and when $k\geq 2$ it gives
\begin{equation}\label{eqs0311b}
 \begin{aligned}
    |\partial_t^2 u_{k\tau}(x,t)|
\leq& C\left[1+\left(t-\tau\right)^{\alpha-1}+\int_{0}^{\tau}\|\partial_sF_{k\tau}(x,s)\|_{G^{1/2}}ds+\int_{\tau}^{t}\|\partial_sF_{k\tau}(x,s)\|_{G^{1/2}}ds\right]\\
\leq& C\left[1+(t-\tau)^{\alpha-1}\right],\\
\end{aligned}
\end{equation}
where the fact that $k\geq 2$ means $t>\tau$ is used.

Next we discuss the boundedness of $\partial_{x}^{n}u_{k\tau}(x,t)$, $n=1,2.$ By \cite[Pages 151--152]{Luchko2012Initial}, it is easy to know that
\begin{equation}\label{eqs2_1209d}
    \begin{aligned}
      \left| \sum_{i=1}^{\infty}\varphi_{i}(0)E_{\alpha}(-\lambda_{i}t^{\alpha})X_{i}^{(n)}(x) \right|\leq C,\  n=1,2.
    \end{aligned}
\end{equation}
In addition, by performing analogous manipulation of \cite[Pages 151--152]{Luchko2012Initial}, it is obvious that
\begin{equation}\label{eqs2_1209e}
    \begin{aligned}
\sum_{i=1}^{\infty}\left|F_{k\tau,i}(t-s)E_{\alpha}(-\lambda_{i}s^{\alpha})
X_{i}^{(n)}(x)\right|\leq C, \ n=1,2.
    \end{aligned}
\end{equation}
Since
\begin{equation}
\begin{aligned}
|\partial_x^n u_{k\tau}(x,t)|&=\left|\sum_{i=1}^{\infty}\varphi_{i}(0)E_{\alpha}(-\lambda_{i}t^{\alpha})X_{i}^{(n)}(x)+\int_{0}^{t}\sum_{i=1}^{\infty}F_{k\tau,i}(t-s)E_{\alpha}(-\lambda_{i}s^{\alpha})X_{i}^{(n)}(x)ds\right|\\&\leq \left|\sum_{i=1}^{\infty}\varphi_{i}(0)E_{\alpha}(-\lambda_{i}t^{\alpha})X_{i}^{(n)}(x)\right|+\int_{0}^{t}\left|\sum_{i=1}^{\infty}F_{k\tau,i}(t-s)E_{\alpha}(-\lambda_{i}s^{\alpha})X_{i}^{(n)}(x)\right|ds,
\end{aligned}
\end{equation}
we can immediately obtain
$|\partial_x^n u_{k\tau}(x,t)|\leq C,\ n=1, 2$ from (\ref{eqs2_1209d}) and (\ref{eqs2_1209e}).

Now we discuss (\ref{eqs4_1}) and (\ref{eqs4_2}).
For $t\in((k-1)\tau,k\tau]$, let $$Y_1(x,t)=\sum_{i=1}^\infty\int_0^tE_\alpha(-\lambda_is^\alpha)F_{k\tau,i}(t-s)ds\cdot X_i-E_\alpha(0)F_\tau(x,0)t, Y_2(x,t)=\sum_{i=1}^\infty\varphi_i(0)\sum_{j=m+1}^\infty\frac{(-\lambda_it^\alpha)^j}{\Gamma(j\alpha+1)}\cdot X_i.$$
Then it follows from (\ref{eqs2_5}) and (\ref{eqs2_6}) that the solution of (\ref{eqs1_1})--(\ref{eqs1_3}) can be written as
\begin{equation}\label{eqs2_250323a}
\begin{aligned}
u_{k\tau} (x,t)=&\sum_{j=0}^{m}\frac{(-1)^jt^{j\alpha}}{\Gamma(j\alpha+1)}G^j\left(\varphi(x,0)\right)+E_\alpha(0)F_\tau(x,0)t+Y(x,t),\\
     :=& \sum_{j=0}^{m}\gamma_jt^{j\alpha}+\tilde{\gamma}t+Y(x,t),\ t\in((k-1)\tau,k\tau],
\end{aligned}
\end{equation}
where $Y(x,t)=Y_1(x,t)+Y_2(x,t)$. First, it is obvious that $\gamma_j$ is bounded, and the boundedness of Mittag-Leffler function and the continuity of $\varphi,f$ can ensure the boundedness of $\tilde{\gamma}$.
According to the definition of $Y_1(x,t)$ and $Y_2(x,t)$, it is easy to check $Y(\cdot,0)=0$.
Furthermore, it follows from (\ref{eqs2_5}) and (\ref{eqs2_6}) that
\begin{align*}
    Y_1(x,t)=u_{k\tau}(x,t)-\sum_{i=1}^{\infty}\varphi_{i}(0)E_{\alpha}(-\lambda_{i}t^{\alpha})X_{i}(x)-E_{\alpha}(0)F_\tau(x,0)t,\ (x,t)\in[0,L]\times ((k-1)\tau,k\tau].
\end{align*}
Taking the first derivative on both sides of the above equation on $t$ and according to (\ref{eqs2_12}), it yields
\begin{align*}
    \partial_tY_1(x,t)=&\sum_{i=1}^{\infty}\Big\{E_{\alpha}(0)F_{k\tau,i}(t)+\int_{0}^{t}(-\lambda_i  s^{\alpha-1})E_{\alpha,\alpha}(-\lambda_{i}s^{\alpha})F_{k\tau,i}(t-s)ds\Big\} X_{i}(x)\\&-E_{\alpha}(0)F_\tau(x,0),\ (x,t)\in[0,L]\times ((k-1)\tau,k\tau].
\end{align*}
Besides, it follows from (\ref{eqs2_5}), (\ref{eqs2_6}) and (\ref{eqs2_7}) that $F_{k\tau}$ is obviously a continuous function on $t.$  Therefore
$$
\partial_t Y_1(\cdot,t)|_{t=0}=0,\ Y_1(\cdot,t)\in C^1[0,K\tau] \ \text{and}\left|\partial_tY_1(\cdot,t)\right|\leq C. $$
Note that
\begin{equation}
    \begin{aligned}
 \partial_t^2 Y_1(x,t) =&\sum_{i=1}^{\infty}\Big\{E_{\alpha}(0)\partial_tF_{k\tau,i}(t)
 +(-\lambda_{i}t^{\alpha-1}E_{\alpha,\alpha}(-\lambda_{i}t^{\alpha}))F_{k\tau,i}(0)\\
 &+\int_{0}^{t}-\lambda_i s^{\alpha-1}E_{\alpha,\alpha}(-\lambda_{i}s^{\alpha})
\partial_tF_{k\tau,i}(t-s)ds \Big\}X_i(x)\\
\leq& C\left[\|\partial_tF_{k\tau}(x,t)\|_{G^{1/2}}+\|F_{k\tau}(x,0)\|_{G^{3/2}}t^{\alpha-1}+\int_{0}^{t}\|\partial_tF_{k\tau}(x,t-s)\|_{G^{1/2}}ds \right].
    \end{aligned}
\end{equation}
Hence Lemma \ref{jullemma2.2} implies that
$\partial_t^2 Y_1(x,t) \leq C(1+t^{\alpha-1}),\ t\in(0,\tau]$ and $\partial_t^2 Y_1(x,t) \leq C(1+(t-\tau)^{\alpha-1}),\ t\in(\tau,K\tau].$
For $\partial_tY_2(x,t)$, it is clear that  $$\partial_tY_2(x,t)=\sum_{i=1}^\infty\varphi_i(x,0)(-\lambda_i)^{m+1}t^{(m+1)\alpha-1} X_i\cdot\left(\sum_{j=0}^\infty\frac{(-\lambda_i)^jt^{j\alpha}}{\Gamma(j\alpha+(m+1)\alpha)}\right).$$ The boundedness of the Mittag-Leffler function gives
$$|\partial_tY_2(x,t)|\leq C\left|G^{m+1}(\varphi(x,0))\right|t^{(m+1)\alpha-1} . $$
Similarly, we can also obtain
$$|\partial^2_tY_2(x,t)|\leq C\left|G^{m+1}(\varphi(x,0))
 \right|t^{(m+1)\alpha-2}.$$
Therefore, when $(m+1)\alpha>1$, one has $$\partial_t Y_2(\cdot,t)|_{t=0}=0,Y_2(\cdot,t)\in C^1[0,K\tau] \ \text{and}\ \left|\partial_t^2Y_2(\cdot,t)\right|\leq Ct^{(m+1)\alpha-2}. $$The above discussion means that (\ref{eqs4_1}) and (\ref{eqs4_2}) hold.
$\quad\Box$

\section{Fully discrete scheme and its stability}
\label{Sec3}
Similar to \cite[Eq. (3.1)]{Bu2024Finite} to
take the Riemann-Liouville integral ${}_0I_t^{1-\alpha}$ to the both sides of (\ref{eqs1_1}), then by \cite[Property 2.4]{li2019Theory} the problem (\ref{eqs1_1})--(\ref{eqs1_3}) can be transformed into
\begin{equation}\label{eqs3_1}
\begin{aligned}
{}^C_0D_{t}^{\alpha}u(x,t)=&p\Delta u(x,t)+au(x,t)+b{}_0I_t^{1-\alpha}u(x,t-\tau)+\tilde{f}(x,t),\quad (x,t)\in \Omega \times (0,K\tau],\\
\end{aligned}
\end{equation}
\begin{equation}\label{eqs3_2}
    \quad \quad \quad u(x,t)=\varphi(x,t),\quad (x,t)\in \bar{\Omega} \times \left[-\tau,0\right],
\end{equation}
\begin{equation}\label{eqs3_3}
 u(x,t)|_{\partial \Omega}=0, \quad t \in [-\tau,K\tau],
\end{equation}
where ${}_0I_t^{\beta}u(x,t):=\int_0^{t}\omega_{\beta}(t-s)u(x,s)ds$, $\tilde{f}(x,t)={}_0I_t^{1-\alpha}f(x,t)$ and ${}^C_0D_{t}^{\alpha}u(x,t)$ is the Caputo fractional derivative defined by $${}^C_0D_{t}^{\alpha}u(x,t)=\int_0^{t}\omega_{1-\alpha}(t-s)\partial_su(x,s)ds.$$
In order to establish a suitable numerical scheme to facilitate the discussion of error estimates, when $t> \tau,$ we rewrite the third term of the right hand of (\ref{eqs3_1}) as
\begin{equation}
    \begin{aligned}
        {}_0I_t^{1-\alpha}u(x,t-\tau)=&\int_\tau^{t}\omega_{1-\alpha}(t-s)u(x,s-\tau)ds+\int_0^{\tau}\omega_{1-\alpha}(t-s)\varphi(x,s-\tau)ds\\=&\int_0^{t-\tau}\omega_{1-\alpha}(t-\tau-s)u(x,s)ds+\int_0^{t}\omega_{1-\alpha}(t-s)\tilde{\varphi}(x,s-\tau)ds\\&-\partial_t\varphi(x,0)\int_\tau^{t}\omega_{1-\alpha}(t-s)\left(s-\tau\right)ds-\varphi(x,0)\int_\tau^{t}\omega_{1-\alpha}(t-s)ds\\=&\int_0^{t-\tau}\omega_{1-\alpha}(t-\tau-s)\left(u(x,s)-\varphi(x,0)\right)ds+\int_0^{t}\omega_{1-\alpha}(t-s)\left(\tilde{\varphi}(x,s-\tau)-\varphi(x,-\tau)\right)ds\\&-\left[\frac{(t-\tau)^{2-\alpha}}{\Gamma(3-\alpha)}\partial_t\varphi(x,0)-\frac{t^{1-\alpha}}{\Gamma(2-\alpha)}\varphi(x,-\tau)\right],\\
    \end{aligned}
\end{equation}
where
\begin{equation}
    \begin{aligned}
     \tilde{\varphi}(x,t-\tau)=\left\{ \begin{aligned}
		&\varphi(x,t-\tau), \ &t\leq \tau,\\
        &\partial_t\varphi(x,0)(t-\tau)+\varphi(x,0),\ &t> \tau.\\
	\end{aligned} \right.
    \end{aligned}
\end{equation}
The above equality implies that
\begin{equation}
\begin{aligned}
     {}_0I_t^{1-\alpha}u(x,t-\tau)=\left\{ \begin{aligned}
	 &{}_0I_t^{1-\alpha}\left(\varphi(x,t-\tau)-\varphi(x,-\tau)\right)+\frac{t^{1-\alpha}}{\Gamma(2-\alpha)}\varphi(x,-\tau),\ &t\leq \tau,\\	
        &{}_0I_{t-\tau}^{1-\alpha}\left(u(x,t)-\varphi(x,0)\right)+{}_0I_t^{1-\alpha}\left(\tilde{\varphi}(x,t-\tau)-\varphi(x,-\tau)\right)\\&-\left[\frac{(t-\tau)^{2-\alpha}}{\Gamma(3-\alpha)}\partial_t\varphi(x,0)-\frac{t^{1-\alpha}}{\Gamma(2-\alpha)}\varphi(x,-\tau)\right], \ &t> \tau.\\
	\end{aligned} \right.&
\end{aligned}
\end{equation}
Therefore, for $t=t_n,$ (\ref{eqs3_1}) can be rewritten as
\begin{equation}\label{eqs106a}
\begin{aligned}
{}^C_0D_{t}^{\alpha}u(x,t_n)=&p\Delta u(x,t_n)+au(x,t_n)+b{}_0\bar{I}_t^{1-\alpha}u(x,t_n-\tau)+F^n(x),
\end{aligned}
\end{equation}
where $${}_0\bar{I}_t^{1-\alpha}u(x,t_n-\tau)=\left\{ \begin{aligned}
		&{}_0I_{t}^{1-\alpha}\left(\varphi(x,t_n-\tau)-\varphi(x,-\tau)\right), \ &t_n\leq \tau,\\
        &{}_0I_t^{1-\alpha}\left(\tilde{\varphi}(x,t_n-\tau)-\varphi(x,-\tau)\right)+{}_0I_{t-\tau}^{1-\alpha}\left(u(x,t_n)-\varphi(x,0)\right),\ &t_n> \tau,\\
	\end{aligned} \right.
    $$ and $$F^n(x)=\left\{ \begin{aligned}
		&\tilde{f}^n(x)+b\frac{t_n^{1-\alpha}}{\Gamma(2-\alpha)}\varphi(x,-\tau), \ &t_n\leq \tau,\\
        &\tilde{f}^n(x)-b\left[\frac{(t_n-\tau)^{2-\alpha}}{\Gamma(3-\alpha)}\partial_t\varphi(x,0)-\frac{t_n^{1-\alpha}}{\Gamma(2-\alpha)}\varphi(x,-\tau)\right],\ &t_n> \tau.\\
	\end{aligned} \right.$$

Now we consider the discretization of fractional operators. For a given
integer $N>0$, let $\rho=\frac{\tau}{N}$, $t_n=n\rho,\ -N\leq n\leq KN,$ $u^k=u(t_k)$ and $\nabla_tu^k=u^k-u^{k-1},\ 1\leq
k\leq KN,$ where $u(t_k)$ denotes the value of $u$ at $t=t_k$.
According to \cite[Remark 2.13]{Kilbas2006Theory} and Theorem \ref{thm2.1}, we notice that
\begin{equation}\label{eqs0124f}
{}^C_0D_{t}^{\alpha}u(t)=\partial_t^{\alpha}(u(t)-u(0)).
\end{equation}
Therefore, the Gr\"unwald-Letnikov type approximation \cite[Section 5.3]{li2019Theory} gives \begin{equation}\label{eqs3_5}
    \begin{aligned}
{}^C_0D_{t}^{\alpha}u^n&=\rho ^{-\alpha}\sum_{k=0}^ng_k^{(\alpha)}(u^{n-k}-u^0)+r_1^n,
    \end{aligned}
\end{equation}
where $g_k^{(\alpha)}:=(-1)^k \left(\begin{array}{c}{\alpha}\\{k}\end{array}\right)$ and $r_1^n$ is the corresponding local truncation error. Furthermore, we generalize the above Gr\"unwald-Letnikov type approximation to discretize the fractional integral operators. If $t_n\leq\tau$, then
\begin{equation}\label{eqs122}
    {}_0I_{t}^{1-\alpha}\left(\varphi(t_n-\tau)-\varphi(-\tau)\right)=\rho ^{1-\alpha}\sum _{k=0}^{n}g_{n-k}^{(\alpha-1)}\left(\varphi(t_k-\tau)-\varphi(-\tau)\right)+r_2^{n-N};
\end{equation}
and if $t_n> \tau,$ then
\begin{equation}\label{eqs3_6}
\begin{aligned}
&{}_0I_t^{1-\alpha}\left(\tilde{\varphi}(t_n-\tau)-\varphi(-\tau)\right)+{}_0I_{t-\tau}^{1-\alpha}\left(u(t_n)-\varphi(0)\right)\\=&\rho ^{1-\alpha}\sum _{k=0}^{n}g_{n-k}^{(\alpha-1)}\left(\tilde{\varphi}(t_k-\tau)-\varphi(-\tau)\right)+\rho ^{1-\alpha}\sum _{k=0}^{n-N}g_{n-N-k}^{(\alpha-1)}\left(u(t_k)-\varphi(0)\right)+r_2^{n-N}\\=&\rho ^{1-\alpha}\sum _{k=0}^{N-1}g_{n-k}^{(\alpha-1)}\tilde{\varphi}(t_k-\tau)+\rho ^{1-\alpha}\sum _{k=0}^{n-N}g_{n-N-k}^{(\alpha-1)}u(t_k)+\rho ^{1-\alpha}\sum _{k=N}^{n}g_{n-k}^{(\alpha-1)}\tilde{\varphi}(t_k-\tau)\\&-\varphi(-\tau)\cdot\rho ^{1-\alpha}\sum _{k=0}^{n}g_{n-k}^{(\alpha-1)}-\varphi(0)\cdot\rho ^{1-\alpha}\sum _{k=0}^{n-N}g_{n-N-k}^{(\alpha-1)}+r_2^{n-N}\\=&\rho^{1-\alpha}\sum_{k=0}^ng_{n-k}^{(\alpha-1)}u^{k-N}+\partial_t\varphi(0)\cdot\rho^{1-\alpha}\sum _{k=N}^{n}g_{n-k}^{(\alpha-1)}(t_k-\tau)-\varphi(-\tau)\cdot\rho ^{1-\alpha}\sum _{k=0}^{n}g_{n-k}^{(\alpha-1)}+r_2^{n-N}.
    \end{aligned}
\end{equation}

It should be pointed out that the above discretization of fractional integrals can be considered as a special case of Lubich's convolution quadrature method \cite{Lubich1986Discretized,Lubich1988Convolution}, and
(\ref{eqs122}) and (\ref{eqs3_6}) can be uniformly written as
\begin{equation}\label{eqs124a}
    {}_0\bar{I}_t^{1-\alpha}u^{n-N}=J^{1-\alpha}u^{n-N}+H^n+r_2^{n-N},
\end{equation}
where $J^{1-\alpha}u^{n-N}=\rho^{1-\alpha}\sum_{k=0}^ng_{n-k}^{(\alpha-1)}u^{k-N}$, $H^n=\partial_t\varphi(0)\cdot\rho^{1-\alpha}\sum _{k=N}^{n}g_{n-k}^{(\alpha-1)}(t_k-\tau)-\varphi(-\tau)\cdot\rho ^{1-\alpha}\sum _{k=0}^{n}g_{n-k}^{(\alpha-1)}$ and the local truncation error $r_2^{n-N}={}_0\bar{I}_t^{1-\alpha}u^{n-N}-J^{1-\alpha}u^{n-N}-H^n.$ In addition, according to \cite{Chen2022Sharp}, one has $g_k^{(\alpha)}=\rho ^{\alpha}\left(A_k-A_{k-1}\right)$ with $A_{k-1}=\rho^{-\alpha}\frac{\Gamma(k-\alpha)}{\Gamma(1-\alpha) \Gamma(k)}
,k\geq 1.$
Hence (\ref{eqs3_5}) can be further rewritten as
\begin{equation}\label{eqs124b}
    \begin{aligned}
{}^C_0D_{t}^{\alpha}u(t_n)&=\sum_{k=1}^nA_{n-k}\nabla_tu^k+r_1^n\\&:=\bar{\partial}_\rho ^{\alpha}u^n+r_1^n.
    \end{aligned}
\end{equation}
We also note that $g_k^{(\alpha-1)}=\frac{\Gamma(k+1-\alpha)}{\Gamma(1-\alpha)\Gamma(k+1)}.$
It means that
\begin{equation}
    J^{1-\alpha}u^{n-N}=\rho\sum_{k=0}^nA_ku^{n-N-k}.
\end{equation}
Let $X_h$ be a continuous piecewise linear finite element space of $H_0^1(\Omega)$ under the quasi-uniform partition on $\Omega$ with the maximum diameter $h$. By (\ref{eqs124a}) and (\ref{eqs124b}), we propose the fully discrete finite element scheme of problem (\ref{eqs3_1})--(\ref{eqs3_3}): find $u^n_h\in X_h$, $n=1,2,\cdots,KN$ such that \begin{equation}\label{eqs3_9}
\left(\bar{\partial}_\rho^{\alpha} u_h^{n},v_h\right)=-p(\nabla u_h^n,\nabla v_h)+a(u_h^n,v_h)+b\left(J^{1-\alpha}u_h^{n-N},v_h\right)+b\left(H^n,v_h\right)+\left(F^n,v_h\right),\ \forall v_h\in X_h,
\end{equation}
where $u_h^n\in X_h, -N\leq n\leq 0$ is a suitable approximation of the initial value function $\varphi(t_n)$.

In order to discuss the stability of the fully discrete finite element scheme (\ref{eqs3_9}), we introduce three useful lemmas firstly.
\begin{lem}\label{juellemma3_2}
For the coefficients $A_n,\ 1\leq n\leq KN$, one has
$$A_n>0, \ \rho\sum_{k=0}^nA_k<\rho^{1-\alpha}+\frac{t_n^{1-\alpha}}{\Gamma(2-\alpha)}.$$
\end{lem}
\noindent\textbf{Proof.}
According to the definition of $A_n,$ it is obvious that $A_n>0$ and
\begin{equation}
\rho\sum_{k=0}^nA_k=\rho^{1-\alpha}\sum_{k=1}^{n+1}\frac{\Gamma(k-\alpha)}{\Gamma(1-\alpha)\Gamma(k)}.
\end{equation}
From \cite[(3) and (4)]{Chen2022Sharp}, we know that
\begin{equation}
    \frac{\Gamma(k-\alpha)}{\Gamma(1-\alpha)\Gamma(k)}<\frac{(k-1)^{-\alpha}}{\Gamma(1-\alpha)}\ \text{for}\ k\geq 2.
\end{equation}
Therefore
 \begin{equation}
    \begin{aligned}
         \rho\sum_{k=0}^nA_k&<\rho^{1-\alpha}+\rho^{1-\alpha}\sum_{k=2}^{n+1}\frac{(k-1)^{-\alpha}}{\Gamma(1-\alpha)}\\&<\rho^{1-\alpha}+\frac{\rho^{1-\alpha}}{\Gamma(1-\alpha)}\sum_{k=2}^{n+1}\int_{k-2}^{k-1}s^{-\alpha}ds\\&<\rho^{1-\alpha}+\frac{t_n^{1-\alpha}}{\Gamma(2-\alpha)}.
   \quad\Box \end{aligned}
 \end{equation}
Define
\begin{equation}\label{eqs3_10}       P_0=\frac{1}{A_0},\ P_{n}=\frac{1}{A_0}\sum_{j=0}^{n-1}P_{j}\left(A_{n-j-1}-A_{n-j}\right),
\end{equation}
and
\begin{equation}
    \begin{aligned}
     K_{\beta,n}=\left\{ \begin{aligned}
		&1+\frac{1-n^{1-\beta}}{\beta-1}, &\beta\neq 1,\\
        &1+\ln n ,\ \ \ & \beta=1.\\
	\end{aligned} \right.
    \end{aligned}
\end{equation}
Then the following inequality and properties hold.
\begin{lem}\label{juellemma3_3}\cite{Chen2022Sharp,Liao2019A}
Let $\left\{g^j\right\}_{j=1}^{KN}$ and $\left\{\mu_j\right\}_{j=0}^{KN}$ be given nonnegative sequences. Suppose that there exists a constant $\Lambda$ such that $\Lambda\geq \sum_{j=0}^{KN}\mu_j$ and the step size satisfies $\rho\cdot\sqrt[\alpha]{2\pi_A\Gamma(2-\alpha)\Lambda}\leq 1.$
Then, for any nonnegative sequence $\left\{v^j\right\}_{ j=0}^{KN}$ such that
 \begin{equation}
      \sum_{k=1}^nA_{n-k}\nabla_tv^k\leq  \sum_{k=1}^n\mu_{n-k}v^k+g^n \ \text{for} \ 1\leq n\leq KN,
 \end{equation}
 it holds that
 \begin{equation}
     v^n\leq 2E_\alpha\left(2\pi_A\Lambda t_n^{\alpha}\right)\left(v^0+\max_{1\leq k\leq n}\sum_{j=1}^kP_{k-j}g^j\right)\ \text{for} \ 1\leq n\leq KN,
 \end{equation}
where the value of $\pi_A,$ i.e. $\pi_A=2^{\alpha}$ has been discussed in \cite{Chen2022Sharp}.
\end{lem}
\begin{lem}\label{juellemma3_1}\cite{Chen2022Sharp}
For the discrete coefficient $P_j,$ it has the following properties\\
\begin{align*}
  &(\romannumeral1)  \ \sum_{j=1}^nP_{n-j}\leq \frac{2^{\alpha}t_n^{\alpha}}{\Gamma(1+\alpha)};\\
  &(\romannumeral2)  \ \sum_{j=k}^nP_{n-j}A_{j-k}=1,\ 1\leq k\leq n;\\
&(\romannumeral3)  \ \sum_{j=1}^nj^{-\beta}P_{n-j}\leq \rho ^\alpha n^{-\beta}+\frac{\rho ^\alpha}{\Gamma(\alpha)}\left[K_{\beta,n}\left(\frac{n}{2}\right)^{\alpha-1}+\frac{1}{\alpha}\left(\frac{n}{2}\right)^{\alpha-\beta}\right],
\end{align*}
where $\beta\geq 0$ is a constant and $n=1,2,\cdots,KN.$
\end{lem}

Now we state the stability of the fully discrete finite element scheme (\ref{eqs3_9}).
\begin{thm}
Let $u_h^n$ be the solution to (\ref{eqs3_9}), and denote $c=\sup_{0\leq k\leq N}\left\|u_h^{k-N}\right\|_0$. Then for $\ 1\leq n\leq KN$, it holds
\begin{equation}\label{eqs0124e}
    \begin{aligned}
\left\|u_h^n\right\|_0\leq 2E_{\alpha}\left(2^{\alpha+1}\Lambda t_n^\alpha\right)\left(\left\|u_h^0\right\|_0+\frac{2^{\alpha}t_n^{\alpha}}{\Gamma(1+\alpha)}\max_{1\leq k\leq n}y^k\right) ,
    \end{aligned}
\end{equation}
 where if $1\leq n\leq N,$ then
$$\Lambda=0,\ y^k=|b|\left(\rho^{1-\alpha}+\frac{t_k^{1-\alpha}}{\Gamma(2-\alpha)}\right)\left(c+\left\|\varphi(-\tau)\right\|_0\right)+\left\|F^k\right\|_0;$$
and if $N+1\leq n\leq KN,$ then
$$\Lambda=|b|\left(\rho^{1-\alpha}+\frac{(K\tau)^{1-\alpha}}{\Gamma(2-\alpha)}\right),\ y^k= c|b|\rho \sum_{j=0}^NA_{k-j}+|b|\left\|H^k\right\|_0+\left\|F^k\right\|_0.$$
\end{thm}
\noindent\textbf{Proof.}
Taking $v_h=u_h^n$ in (\ref{eqs3_9}), then
\begin{equation}
\left(\bar{\partial}_\rho^{\alpha} u_h^{n},u_h^n\right)=-p(\nabla u_h^n,\nabla u_h^n)+a(u_h^n,u_h^n)+b\left(J^{1-\alpha}u_h^{n-N},u_h^n\right)+b\left(H^n,u_h^n\right)+\left(F^n,u_h^n\right),
\end{equation}
i.e.
\begin{equation}\label{eqs0124b}
\begin{aligned}
A_0\left\|u_h^n\right\|^2_0\leq& \sum_{k=1}^{n-1}(A_{k-1}-A_k)\left\|u_h^{n-k}\right\|_0\left\|u_h^n\right\|_0+A_{n-1}\left\|u_h^0\right\|_0\left\|u_h^n\right\|_0-p\left\|\nabla u_h^n\right\|^2_0+a\left\|u_h^n\right\|_0^2\\&+|b|\left\|J^{1-\alpha}u_h^{n-N}\right\|_0\left\|u_h^n\right\|_0+|b|\left\|H^n\right\|_0\left\|u_h^n\right\|_0+\left\|F^n\right\|_0\left\|u_h^n\right\|_0.
\end{aligned}
\end{equation}
Using the fact $p>0,$ $a\leq 0$ and $A_{k-1}>A_k$ mentioned in \cite{Chen2022Sharp}, (\ref{eqs0124b}) leads to
\begin{equation}\label{eqs107a}
    \begin{aligned}
        \sum_{k=1}^nA_{n-k}\nabla_t\left\|u_h^k\right\|_0\leq&|b|\left\|J^{1-\alpha}u_h^{n-N}\right\|_0+|b|\left\|H^n\right\|_0+\left\|F^n\right\|_0.
    \end{aligned}
\end{equation}
It follows from Lemma \ref{juellemma3_2} and (\ref{eqs107a}) that if $1\leq n\leq N,$ then
\begin{equation}\label{eqs0124c}
    \begin{aligned}
        \sum_{k=1}^nA_{n-k}\nabla_t\left\|u_h^k\right\|_0&\leq|b|\rho\sum_{k=0}^nA_k\left\|u_h^{n-N-k}\right\|_0+|b|\rho\sum_{k=0}^nA_k\left\|\varphi(-\tau)\right\|_0+\left\|F^n\right\|_0\\&\leq |b|\left(\rho^{1-\alpha}+\frac{t_n^{1-\alpha}}{\Gamma(2-\alpha)}\right)\left(c+\left\|\varphi(-\tau)\right\|_0\right)+\left\|F^n\right\|_0;
    \end{aligned}
\end{equation}
and if $ N+1\leq n\leq KN,$ then
\begin{equation}\label{eqs0124d}
    \begin{aligned}
        \sum_{k=0}^nA_{n-k}\nabla_t\left\|u_h^k\right\|_0&\leq|b|\rho \sum_{k=0}^NA_{n-k}\left\|u_h^{k-N}\right\|_0+|b|\rho\sum_{k=N+1}^nA_{n-k}\left\|u_h^{k-N}\right\|_0+|b|\left\|H^n\right\|_0+\left\|F^n\right\|_0\\&\leq |b|\rho \sum_{k=1}^{n-N}A_{n-N-k}\left\|u_h^{k}\right\|_0+\left(c|b|\rho \sum_{k=0}^NA_{n-k}+|b|\left\|H^n\right\|_0+\left\|F^n\right\|_0\right).
    \end{aligned}
\end{equation}
Applying Lemma \ref{juellemma3_3} to (\ref{eqs0124c}) and (\ref{eqs0124d}), and
 noting that inequalities mentioned in Lemmas \ref{juellemma3_2} and \ref{juellemma3_1}, it yields (\ref{eqs0124e}) immediately.
$\quad\Box$

\section{Convergence analysis}
\label{Sec4}
In this section, we consider the convergence of the fully discrete scheme (\ref{eqs3_9}). Now we give three useful lemmas.
\begin{lem}\label{juellemma5_1}\cite[Lemma 2.1]{Zeng2017Second}
Let $w(t)=t^\sigma,$ $\sigma>0$ and $\beta$ is a real number. Then
\begin{equation}
\left|{}_{RL}D_{0,t}^{\beta}w^n-\rho^{-\beta}\sum_{j=0}^{n}g_j^{(\beta)}w^{n-j}\right|\leq C\left(\frac{\beta}{2}\left|\frac{\Gamma(\sigma+1)}{\Gamma(\sigma-\beta)}\right|\rho t_n^{\sigma-1-\beta}+\rho^2 t_n^{\sigma-2-\beta}\right),
\end{equation}
where ${}_{RL}D_{0,t}^{\beta}w$ represents the Riemann-Liouville fractional derivative $\partial_t^{\beta}w$ for $\beta>0,$ and ${}_{RL}D_{0,t}^{\beta}w$ is the Riemann-Liouville fractional integral ${}_0I_t^{-\beta}w$ for $\beta<0.$
\end{lem}
\begin{lem}\label{juellemma5_2}\cite[Theorem 1]{Zhao2015A} Assume that $w(t)\in C^1[0,K\tau]$, $\partial^2_tw(t)\in L^1[0,K\tau]$ and $w(0)=\partial_tw(0)=0$. Then
\begin{equation}
\left|{}^C_0D_{t}^{\alpha}w^n-\bar{\partial}_\rho ^{\alpha}w^n\right|\leq C\rho.
\end{equation}
\end{lem}

Combining (\ref{eqs4_1}), (\ref{eqs4_2}), (\ref{eqs0124f}),  Lemma \ref{juellemma5_1} and Lemma \ref{juellemma5_2}, it is clear that the local truncation error  $r_1^n$ in (\ref{eqs124b}) can be bounded by
\begin{equation}\label{eqs_1223c}
\begin{aligned}
   |r_1^n|&\leq C
\left[\sum_{j=1}^m\left(\left|\frac{\Gamma(j\alpha+1)}{\Gamma((j-1)\alpha)}\right|\rho t_n^{(j-1)\alpha-1}+\rho^2t_n^{(j-1)\alpha-2} \right)+\rho t_n^{-\alpha}+\rho^2 t_n^{-\alpha-1}+\rho\right]\\&\leq C\left(\rho t_n^{-1}+\rho^2 t_n^{-2}+\rho t_n^{-\alpha} \right)\\&\leq C\rho t_n^{-1}.
\end{aligned}
\end{equation}
\begin{lem}\cite[Theorem 3.1]{Lubich1988Convolution}\label{juellemma5_3}
      Assume that $w(t)\in C^1[0,K\tau],$ then
\begin{equation}
  \left|{}_0I_t^{1-\alpha}w^n-\rho ^{1-\alpha}\sum _{j=0}^{n}g_{n-j}^{(\alpha-1)}w^j\right|\leq Ct_n^{-\alpha}\rho\left(|w^0|+t_n\max_{0\leq s\leq t_n}|\partial_tw(s)|\right).
\end{equation}
\end{lem}

Suppose that
\begin{equation}\label{eqs3_250323b}
\varphi(t)\in C^1([-\tau,0])\ \text{and}\ |\partial_t\varphi(t)|\leq C,
\end{equation}
for all $x\in[0,L].$
Then it follows from Lemma \ref{juellemma5_1} and Lemma \ref{juellemma5_3} that
\begin{equation}\label{eqs106b}
  |r_2^{n-N}|\leq C \rho ,
  \ n\geq 1,
\end{equation}
where the splitting $
    u=\left(\sum_{j=1}^m \gamma_jt^{j\alpha}+\tilde{\gamma}t\right)+\left(\varphi(0)+Y\right)
   $ mentioned in (\ref{eqs4_1}) is used. Furthermore, for $P_j$ defined in (\ref{eqs3_10}), the combination of ($\romannumeral1$) in Lemma \ref{juellemma3_1} and (\ref{eqs106b}) gives
   \begin{equation}
       \sum_{j=1}^{n}P_{n-j}|r_2^{j-N}|\leq  C\rho.
   \end{equation}

In order to obtain a new Gronwall inequality to discuss the convergence of the numerical scheme (\ref{eqs3_9}), we give the following lemma.
\begin{lem}\label{jullemma5.6} Let  $Z_1=p_1\left( 1,1,\cdots ,1 \right)^T\in \mathbb{R}^n$, $Z_2=p_2\rho \left( (1+\ln n)t_n^{\alpha-1},(1+\ln (n-1))t_{n-1}^{\alpha-1},\cdots ,t_1^{\alpha-1}\right)^T $ and
\begin{equation}\label{Eq:matrix1}
M=p_3\rho\begin{array}{c}
\left(
\begin{array}{ccccccc}
	0 & \cdots & 0 & 1 & 1 & \cdots & 1 \\[6pt]
	0 & \cdots & 0 & 0 & 1 & \cdots & 1 \\[6pt]
\vdots&\vdots&\vdots&\vdots&\vdots&\ddots&\vdots\\[6pt]
	0 & \cdots & 0 & 0 & 0 & \cdots & 1 \\[6pt]
    0 & \cdots & 0 & 0 & 0 & \cdots & 0 \\[6pt]
			\vdots & \vdots & \vdots & \vdots & \vdots & \vdots & \vdots \\[6pt]
			0 & \cdots & 0 & 0 & 0 & \cdots & 0
		\end{array}
		\right)_{n \times n}
		\\[-1.0em]
		\begin{array}{ccccccc}
			\multicolumn{3}{c}{\hspace{-1.7em}\underbrace{\phantom{0 \quad \cdots \quad 0}}_{N}}
			& \phantom{0} & \phantom{0} & \phantom{\cdots} & \phantom{0}
		\end{array}
	\end{array} ,
\end{equation}
where $p_1,\ p_2,\ p_3$ are some nonnegative constants and $n\in [(k-1)N+1,kN], k\geq 1.$ Then we have\\
$(\romannumeral1)$
$
M^j=\mathbb{O},\quad j\geq k;
$\\
$(\romannumeral2)$ For $ M^q Z_1,$  $q=0,1,\cdots,k-1,$ the following inequalities hold
 \begin{equation}
     M^q Z_1\leq Cp_1\left(1,1,\cdots,1,\underset{qN}{\underbrace{0,\cdots, 0}}\right)^T,
 \end{equation}
and
\begin{equation}\label{eqs4_8}
     \sum_{j=0}^{k-1} M^j Z_1\leq Cp  _1\left(\underset{n-(k-1)N}{\underbrace{k,\cdots k}},\underset{N}{\underbrace{k-1,\cdots k-1}},\cdots,\underset{N}{\underbrace{1,\cdots, 1}}\right)^T;
 \end{equation}
 $(\romannumeral3)$ For $M^q Z_2,$ $q=1,2,\cdots,k-1,$ the following inequalities hold
\begin{equation}\label{eqs4_9}
     M^q Z_2\leq C \rho \left(1+\ln (n-qN) ,1+\ln (n-qN-1),\cdots ,1,\underset{qN}{\underbrace{0,\cdots ,0}}\right)^T,
 \end{equation}
 and
\begin{equation}\label{eqs4_10}
 \begin{aligned}
  \sum_{j=1}^{k-1} M^j Z_2\leq  C\rho\left(\underset{n-(k-1)N}{\underbrace{(k-1)+\sum_{j=1}^{k-1}\ln (n-j N),\cdots,(k-1)+\sum_{j=1}^{k-1}\ln ((k-1-j)N+1)}},\cdots,\underset{N}{\underbrace{1+\ln (2N-1),\cdots,1}},\underset{N}{\underbrace{0,\cdots, 0}}\right)^T.
 \end{aligned}
 \end{equation}
\end{lem}
\noindent\textbf{Proof.}
First, it follows from the mathematical induction that $(\romannumeral1)$ can be obtained easily. For $M^qZ_1,$ we notice that $(\romannumeral2)$ is  obviously true for $q=0.$ Suppose that $(\romannumeral2)$ holds for $q\leq l$. Then
\begin{equation}
    \begin{aligned}
        M^{l+1}Z_1&=M\cdot M^lZ_1\\
        &\leq Cp _1M\left(1,1,\cdots,1,\underset{lN}{\underbrace{0,\cdots, 0}}\right)^T\\
        &\leq Cp_1\rho \left(n-(l+1)N,n-(l+1)N-1,\cdots,1,\underset{(l+1)N}{\underbrace{0,\cdots 0}}\right)^T\\&\leq Cp  _1\left(1,1,\cdots,1,\underset{(l+1)N}{\underbrace{0,\cdots, 0}}\right)^T,\\
    \end{aligned}
\end{equation}
and
\begin{equation}
     \sum_{j=0}^{k-1} M^j Z_1\leq Cp_1\left(\underset{n-(k-1)N}{\underbrace{k,\cdots k}},\underset{N}{\underbrace{k-1,\cdots k-1}},\cdots,\underset{N}{\underbrace{1,\cdots, 1}}\right)^T.
 \end{equation}
 Therefore the mathematical induction implies that $(\romannumeral2)$ is true.

Now we discuss $(\romannumeral3).$ For $q=1,$ one has
 \begin{equation}
 \begin{aligned}
     MZ_2
&=\lambda p_2 \rho^2 \left(\sum_{j=1}^{n-N} \left(1+\ln j\right)t_j^{\alpha-1},\sum_{j=1}^{n-N-1} \left(1+\ln j\right)t_j^{\alpha-1},\cdots ,t_{1}^{\alpha-1},\underset{N}{\underbrace{0,\cdots ,0}}\right)^T\\
&\leq C \rho^{1+\alpha} \left(\left(1+\ln (n-N)\right)\sum_{j=1}^{n-N} j^{\alpha-1},\left(1+\ln (n-N-1)\right)\sum_{j=1}^{n-N-1} j^{\alpha-1},\cdots ,1,\underset{N}{\underbrace{0,\cdots ,0}}\right)^T.\\
 \end{aligned}
 \end{equation}
Since $\rho^\alpha\sum_{j=1}^{n-N}j^{\alpha-1}\leq \rho^{\alpha}\left(1+\int_{0}^{n-N}s^{\alpha-1}ds\right)\leq C,$ $MZ_2$ can be bounded by
 \begin{equation}
 \begin{aligned}
     MZ_2
&\leq C \rho\left(1+\ln (n-N),1+\ln (n-N-1),\cdots ,1,\underset{N}{\underbrace{0,\cdots ,0}}\right)^T.
 \end{aligned}
 \end{equation}
 Suppose that (\ref{eqs4_9}) holds for $q\leq l.$ When $q=l+1,$ then
 \begin{equation}
 \begin{aligned}
M^{l+1}Z_2&=M\cdot M^{l}Z_2\\
&\leq C\rho^2\left(\sum_{j=1}^{n-(l+1)N} \left(1+\ln j\right),\sum_{j=1}^{n-(l+1)N-1}  \left(1+\ln j\right),\cdots ,1,\underset{(l+1)N}{\underbrace{0,\cdots ,0}}\right)^T\\
&\leq C \rho \left(1+\ln (n-(l+1)N) ,1+\ln (n-(l+1)N-1),\cdots ,1,\underset{(l+1)N}{\underbrace{0,\cdots ,0}}\right)^T.
 \end{aligned}
 \end{equation}
Thus it is clear that $(\romannumeral3)$ is true from the mathematical induction.
$\quad\Box$

Using the above results, we state a new Gronwall inequality as follows.
\begin{lem}\label{juellemma4_4}
     Assume that $\left\{z^n\right\}$ is a non-negative real sequence and satisfies
 \begin{equation}\label{eqs4_16}
        z^n\leq p_3 \rho \sum_{j=N+1}^nz^{j-N}+p_1+p_2\rho\left(1+\ln n\right)t_n^{\alpha-1}, \ (k-1)N+1\leq  n\leq kN,\ k=1,2,\cdots,K,
 \end{equation}
 then
 \begin{equation}
    \begin{aligned}
     z^n\leq C \left[kp_1+\rho \left(\sum_{j=1}^{k-1}\left(1+\ln (n-jN)\right)+p_2\left(1+\ln n\right)t_n^{\alpha-1}\right)\right].
    \end{aligned}
\end{equation}
\end{lem}
\noindent\textbf{Proof.}
For $(k-1)N+1\leq n\leq kN,$ let $Z=\left( z^n,z^{n-1},\cdots ,z^1 \right)^T.$ Then it follows from
 (\ref{eqs4_16}) that
 \begin{equation}
     \begin{aligned}
         Z\leq& MZ+Z_1+Z_2.
     \end{aligned}
 \end{equation}
Applying this inequality repeatedly and noting that the fact $M^k=\mathbb{O},$ it yields
 \begin{equation}
     \begin{aligned}
         Z
         \leq& M\left(MZ+Z_1+Z_2\right)+Z_1+Z_2\\
         \leq& M^2Z+\sum_{j=0}^1M^jZ_1+\sum_{j=0}^1M^jZ_2\\
        &\vdots\\
         \leq& \sum_{j=0}^{k-1}M^jZ_1+\sum_{j=0}^{k-1}M^jZ_2.
     \end{aligned}
 \end{equation}
The combination of (\ref{eqs4_8}) and (\ref{eqs4_10}) gives
 \begin{equation}
     z^n\leq C
\left[kp_1+\rho \left(\sum_{j=1}^{k-1}\left(1+\ln (n-jN)\right)+p_2\left(1+\ln n\right)t_n^{\alpha-1}\right)\right].\quad\Box
\end{equation}

    In order to discuss the spatial error, we introduce the  orthogonal projection operator $P_h:H_0^1(0,L)\rightarrow X _h$ defined by
\begin{equation}
    (\nabla P_h \phi,\nabla v_h)=(\nabla \phi,\nabla v_h),\quad\forall v_h \in X_h.
\end{equation}
In fact, it has a well known property (see \cite[Lemma 1.1]{Thomee2006Galerkin}), i.e.
\begin{equation}\label{eqs_1222b}
  \|P_h \phi-\phi\|_0\leq Ch^2 \|\phi\|_{H^2(0,L)}.
\end{equation}

Now we state the convergence of the numerical scheme (\ref{eqs3_9}).
\begin{thm}\label{th4.1}
Assume that $u(x,t)$ is the solution of (\ref{eqs1_1})--(\ref{eqs1_3}) which can be decomposed into (\ref{eqs4_1})--(\ref{eqs4_2}), $u(x,t)\in H^2(0,L)$ for fixed $t,$ the condition (\ref{eqs3_250323b}) is satisfied, and $u_h^n\in X_h, -N\leq n\leq 0$ in (\ref{eqs3_9}) is a suitable approximation of $\varphi(x,t_n)$ such that $\|u_h^n-\varphi(x,t_n)\|\leq Ch^2.$ Then the numerical solution $u_h^n$ to (\ref{eqs3_9}) satisfies
\begin{equation}
   \left\|u^n-u_h^n\right\|_0\leq C \left[h^2+\rho \left(\sum_{j=1}^{k-1}\left(1+\ln (n-jN)\right)+\left(1+\ln n\right)t_n^{\alpha-1}\right)\right],
\end{equation}
where $(k-1)N+1\leq n\leq kN,\ k=1,2,\cdots,K.$
\end{thm}
\noindent\textbf{Proof.}
First, (\ref{eqs106a}) gives
\begin{equation}\label{eqs4_17}
\begin{aligned}
\left({}^C_0D_{t}^{\alpha}u^n,v_h\right)=&-p\left(\nabla u^n,\nabla v_h\right)+a\left(u^n,v_h\right)+b\left({}_0\bar{I}_t^{1-\alpha}u^{n-N},v_h\right)+\left(F^n,v_h\right).
\end{aligned}
\end{equation}
 Let $\varepsilon_h^n=u_{h}^n-P_h u^n.$ The combination of (\ref{eqs3_9}) and (\ref{eqs4_17}) yields
 \begin{equation}\label{eqs0217b}
 \begin{aligned}
\left(\bar{\partial}_\rho^{\alpha} \varepsilon_h^n,v_h\right)=&-p\left(\nabla \varepsilon_h^n,\nabla v_h\right)+a\left(\varepsilon_h^n,v_h\right)+b\left(J^{1-\alpha}\varepsilon_h^{n-N},v_h\right)+\left(r_1^n,v_h\right)+\left(\bar{\partial}_\rho^{\alpha}(u^n-P_hu^n),v_h\right)\\&-a\left(u^n-P_hu^n,v_h\right)-b\left(r^{n-N}_2,v_h\right)+b\left(J^{1-\alpha}(P_hu^{n-N}-u^{n-N}),v_h\right).
 \end{aligned}
\end{equation}
Taking $v_h=\varepsilon_h^n$ into (\ref{eqs0217b}) leads to
\begin{equation}\label{eqs0217a}
    \begin{aligned}
A_0\|\varepsilon_h^n\|_0^2=&A_{n-1}\left(\varepsilon_h^0,\varepsilon_h^n\right)+\sum_{k=1}^{n-1}\left(A_{n-k-1}-A_{n-k}\right)\left(\varepsilon_h^k,\varepsilon_h^n\right)-p\|\nabla \varepsilon_h^n\|^2_0+a\|\varepsilon_h^n\|^2_0+b\left(J^{1-\alpha}\varepsilon_h^{n-N},\varepsilon_h^n\right)+\left(r_1^n,\varepsilon_h^n\right)\\&+\left(\bar{\partial}_\rho^{\alpha}(u^n-P_hu^n),\varepsilon_h^n\right)-a\left(u^n-P_hu^n,\varepsilon_h^n\right)-b\left(r^{n-N}_2,\varepsilon_h^n\right)+b\left(J^{1-\alpha}(P_hu^{n-N}-u^{n-N}),\varepsilon_h^n\right)\\
\leq &\left[A_{n-1}\left\|\varepsilon_h^0\right\|_0+\sum_{k=1}^{n-1}\left(A_{n-k-1}-A_{n-k}\right)\left\|\varepsilon_h^k\right\|_0+|b|\|J^{1-\alpha}\varepsilon_h^{n-N}\|_0+\|r_1^n\|_0\right.\\&\left.+\left\|\bar{\partial}_\rho^{\alpha}(u^n-P_hu^n)\right\|_0-a\left\|u^n-P_hu^n\right\|_0 +|b|\left\|r^{n-N}_2\right\|_0+|b|\left\|J^{1-\alpha}(P_hu^{n-N}-u^{n-N})\right\|_0\right]\left\|\varepsilon_h^n\right\|_0.
    \end{aligned}
\end{equation}
The above inequality implies that
\begin{equation}\label{eqs0217c}
    \begin{aligned}
        \sum_{i=1}^jA_{j-i}\nabla_t\left\|\varepsilon_h^i\right\|_0\leq& |b|\|J^{1-\alpha}\varepsilon_h^{j-N}\|_0+\|r_1^j\|_0+\left\|\bar{\partial}_\rho^{\alpha}(u^j-P_hu^j)\right\|_0-a\left\|u^j-P_hu^j\right\|_0 +|b|\left\|r^{j-N}_2\right\|_0\\&+|b|\left\|J^{1-\alpha}(P_hu^{j-N}-u^{j-N})\right\|_0.
    \end{aligned}
\end{equation}

Multiplying $P_{n-j}$ on both sides of (\ref{eqs0217c}) and summing on $j$ from $1$ to $n$ give
\begin{equation}\label{eqs4_27}
  \begin{aligned}
\left\|\varepsilon_h^n\right\|_0\leq& |b|\rho \sum_{j=0}^n\left\|\varepsilon_h^{j-N}\right\|_0+\sum_{j=1}^nP_{n-j}\|r_1^j\|_0+\left\|u^n-P_hu^n\right\|_0-a\sum_{j=1}^nP_{n-j}\left\|u^j-P_hu^j\right\|_0\\&+|b|\sum_{j=1}^nP_{n-j}\left\|r^{j-N}_2\right\|_0+|b|\rho \sum_{j=0}^n\left\|P_hu^{j-N}-u^{j-N}\right\|_0+\left\|\varepsilon_h^0\right\|_0.
  \end{aligned}
\end{equation}
For $\sum_{j=1}^nP_{n-j}\|r_1^j\|_0,$ Lemma \ref{juellemma3_1} and (\ref{eqs_1223c}) imply that
\begin{equation}\label{eqs0217e}
\begin{aligned}
    \sum_{j=1}^nP_{n-j}\|r_1^j\|_0&\leq C\sum_{j=1}^nP_{n-j}j^{-1}\\&\leq C\rho^\alpha \left[n^{-1}+\frac{1}{\Gamma(\alpha)}\left((1+\ln n)\left(\frac{n}{2}\right)^{\alpha-1}+\frac{1}{\alpha}\left(\frac{n}{2}\right)^{\alpha-1}\right)\right]
    \\&\leq C\rho \left(1+\ln n\right)t_n^{\alpha-1}.
\end{aligned}
\end{equation}
Similarly we can also easily obtain that
\begin{equation}\label{eqs0217f}
    \sum_{j=1}^nP_{n-j}\|r_2^{j-N}\|_0\leq C\rho.
\end{equation}
Taking the estimates (\ref{eqs0217e}) and  (\ref{eqs0217f}) into (\ref{eqs4_27}), and noting that the fact $\|u^j-u^j_h\|_0\leq Ch^2$ and the assumption $\left\|\varepsilon_h^{j-N}\right\|_0\leq Ch^2$ for $j=0,1\cdots,N,$ one has
\begin{equation}
  \begin{aligned}
\left\|\varepsilon_h^n\right\|_0\leq& |b|\rho \sum_{j=N+1}^n\left\|\varepsilon_h^{j-N}\right\|_0+Ch^2 +C\rho \left(1+\ln n\right)t_n^{\alpha-1}.
  \end{aligned}
\end{equation}
Therefore it follows from Lemma \ref{juellemma4_4} that
\begin{equation}
\left\|\varepsilon_h^n\right\|_0\leq C
\left[h^2+\rho \left(\sum_{j=1}^{k-1}\left(1+\ln (n-jN)\right)+\left(1+\ln n\right)t_n^{\alpha-1}\right)\right].\quad\Box
\end{equation}

In view of Theorem \ref{thm2.1},
we also have an alternative version of the error estimate.
\begin{thm}\label{thm4.2}
Assume that $u(x,t)$ is the solution of (\ref{eqs1_1})--(\ref{eqs1_3}), the conditions of  Theorem \ref{thm2.1} and (\ref{eqs3_250323b}) hold, and $u_h^n\in X_h, -N\leq n\leq 0$ in (\ref{eqs3_9}) is a suitable approximation of $\varphi(x,t_n)$ such that $\|u_h^n-\varphi(x,t_n)\|\leq Ch^2.$ Then the numerical solution $u_h^n$ to (\ref{eqs3_9}) satisfies
\begin{equation}
   \left\|u^n-u_h^n\right\|_0\leq C \left[h^2+\rho \left(\sum_{j=1}^{k-1}\left(1+\ln (n-jN)\right)+\left(1+\ln n\right)t_n^{\alpha-1}\right)\right],
\end{equation}
where $(k-1)N+1\leq n\leq kN,\ k=1,2,\cdots,K.$
\end{thm}
\section{Numerical experiment}
\label{Sec5}
In this section, some numerical tests are presented to support the achieved theoretical results. In the first case, an exact solution is provided for the considered problem, which can
obviously be decomposed
into (\ref{eqs4_1})--(\ref{eqs4_2}). Therefore we will use it to check the Theorem
\ref{th4.1}. Then, in the next case, we do not provide the exact solution, but give
suitable right hand function and initial value function which satisfy the conditions of Theorem \ref{thm4.2}. The obtained numerical results will be
used to test the Theorem \ref{thm4.2}. Define $E(h,N,k)=\max_{(k-1)N+1\leq n\leq kN}\|\bar{u}^n-u_h^n\|_0$ and
$$
rate_{t}:=\log_2\left(\frac{E(h,N)}{E(h,2N)}\right),\  rate_{s}:=\log_2\left(\frac{E(h,N)}{E(h/2,N)}\right),
$$
where $\bar{u}^n=u^n$ when the exact solution is known, otherwise $\bar{u}^n$
represents an approximation of $u^n$ on a sufficiently fine mesh.

\textbf{Example 1.}
For the problem (\ref{eqs1_1})--(\ref{eqs1_3}) on the time and space domain $(0,3]$ and $[0,1]$, two cases are considered as follows:

$\bullet$ In the first case, let $\tau =1 ,$ $p=1/5,$ $a=0,$ $b=1$, and $\varphi(x,t)=-\frac{\Gamma (\alpha+1)}{\pi ^2}sin\pi x$. We take
\begin{align*}
    u(x,t)=\left\{
    \begin{aligned}
        &\left(\frac{-\Gamma(\alpha+1)}{\pi ^2}+t^\alpha\right)sin\pi x,\ &0< t\leq 1,\\
        &\left(\frac{-\Gamma(\alpha+1)}{\pi ^2}+t^\alpha+(t-1)^{\alpha+1}\right)sin\pi x,\ &1< t\leq 2,\\
        &\left(\frac{-\Gamma(\alpha+1)}{\pi ^2}+t^\alpha+(t-1)^{\alpha+1}+(t-2)^{\alpha+2}\right)sin\pi x,\ &2< t\leq 3\\
    \end{aligned}\right.
\end{align*}
as the exact solution. The right hand function $f(x,t)$ can be computed by the exact solution and initial value function.

$\bullet$ In the second case, instead of giving the exact solution, we give
\begin{align*}
    f(x,t)=\left\{
    \begin{aligned}
        &tsin\pi x,\ &0< t\leq 1,\\
        &\left(t+(t-1)^2\right)sin\pi x,\ &1< t\leq 2,\\
        &\left(t+(t-1)^2+(t-2)^{2}\right)sin\pi x,\ &2< t\leq 3,\\
    \end{aligned}\right.
\end{align*}
$\varphi(x,t)=(1+0.1t)\sin{\pi x}$,
and choose the parameters $\tau =1 ,$ $p=1/100,$ $a=-1, b=1$.

The numerical results of the first case are presented in Table \ref{t1} and Table \ref{t2}, and the numerical tests
of the second case are shown in Table \ref{t3} and Table \ref{t4}. In Table \ref{t1} and Table \ref{t3}, we fix  $h=1/1000$ and $h=1/5000$ respectively to ignore the spatial error. Then, $\alpha=0.7$ and $\alpha=0.9$ are chosen for the former,  and $\alpha=0.4$ and $\alpha=0.6$ are taken for the latter to examine the temporal convergence rate. The numerical results of these two tables show that the time convergence order is near $\alpha$ in the first time interval, and the convergence order is $1$ in the following two time intervals. In Table \ref{t2} and Table \ref{t4}, we test the convergence rate at $t=1,2,3$. In these two tables, $h=1/5000$ is chosen to test the convergence in time. It is clear that the results of these two tables are different from those of Table \ref{t1} and Table \ref{t3} because these results show that the convergence order at $t=1,2,3$ is $1$. Then we fix $\rho=1/10000$ and choose different $h$ in Table \ref{t2} and Table \ref{t4}, the results show that the spatial convergence rate is near 2. The above results are obviously agree with the expected cases mentioned in Theorem \ref{th4.1} and Theorem \ref{thm4.2}.

\begin{table}
\caption{Numerical temporal accuracy for different $\alpha.$}\label{t1}
%\resizebox{\textwidth}{40mm}
\begin{center}
\setlength{\tabcolsep}{2.5mm}{
\begin{tabular}{llllllllll}
\hline
$\alpha$& $N$  & $k=1$&   &  $k=2$&  & $k=3$
\\ &  & $E(h,N,k)$ &$rate_t$ & $E(h,N,k)$ &$rate_t$ &  $E(h,N,k)$&$rate_t$ \\
\hline
0.7
&$400$ &9.4632e-04&&4.3779e-04& &1.5550e-03&     \\
&$800$&5.8905e-04& 0.683&2.1786e-04& 1.006&7.7459e-04&1.005    \\
&$1600$ &3.6511e-04&0.690 &1.0774e-04&1.015&3.8416e-04&1.011   \\
&$3200$ &2.2570e-04&0.693&5.2634e-05&1.033&1.8888e-04&1.024    \\
\hline
0.9
&$400$ &1.2244e-04&&6.3351e-04& &2.2168e-03&     \\
&$800$&6.6105e-05&0.889&3.1553e-04&1.005&1.1050e-03&1.004    \\
&$1600$ &3.5548e-05&0.894 &1.5644e-04&1.012&5.4887e-04&1.009   \\
&$3200$ &1.9066e-05&0.898&7.6877e-05&1.024&2.7075e-04&1.019    \\
\hline
\end{tabular}}
\end{center}
\end{table}
%\begin{table}
%\caption{Numerical spatial accuracy for $h=1/M$ and different $\alpha.$}
%%\resizebox{\textwidth}{55mm}
%\begin{center}
%\setlength{\tabcolsep}{2.5mm}{
%\begin{tabular}{llllllllll}
%\hline
%$\alpha$& $M$  & $k=1$&   &  $k=2$&  & $k=3$
%\\ &  & $E(h,N)$ &$rate_t$ & $E(h,N)$ &$rate_t$ &  $E(h,N)$&$rate_t$ \\
%\hline
%0.7&$8$&1.3785e-02&&4.0341e-02&&1.0300e-01&  \\
%&$16$ &3.4689e-03&1.990&1.0150e-02&1.990&2.5921e-02&1.990  \\
%&$32$  &8.6891e-04&1.997&2.5358e-03 &2.000&6.4697e-03&2.002   \\
%&$64$  &2.1759e-04&1.997&6.2783e-04 &2.013&1.5954e-03&2.019  \\
%\hline
%0.9&$8$&1.3269e-02&&4.3445e-02 &&1.1850e-01&  \\
%&$16$ &3.3372e-03&1.991&1.0926e-02 &1.991&2.9812e-02&1.990\\
%&$32$  &8.3567e-04&1.997&2.7271e-03 &2.002&7.4343e-03&2.003   \\
%&$64$  &2.0913e-04&1.998&6.7286e-04 &2.018&1.8271e-03&2.024  \\
%\hline
%\end{tabular}}
%\end{center}
%\end{table}
\begin{table}
\caption{The temporal accuracy and spatial accuracy for $\alpha=0.7$ at $t_{kN}.$}\label{t2}
%\resizebox{\textwidth}{55mm}
\begin{center}
\setlength{\tabcolsep}{2.5mm}{
\begin{tabular}{llllllllll}
\hline
&$k=1$&& $k=2$ &&  $k=3$& \\
\hline
$N$\ \  & $\left\|u^N - u_h^N\right\|_0$ &$rate_t$ & $\left\|u^{2N} - u_h^{2N}\right\|_0$ &$rate_t$ &  $\left\|u^{3N} - u_h^{3N}\right\|_0$&$rate_t$ \\
\hline
400&1.9718e-05&&4.4021e-04&&1.5612e-03&  \\
800&9.7579e-06&1.014&2.2028e-04&0.998&7.8079e-04&0.999  \\
1600&4.8601e-06&1.005&1.1016e-04&0.999&3.9035e-04&1.000  \\
3200&2.4361e-06&0.996&5.5051e-05&1.000&1.9507e-04&1.000  \\
\hline
$h$\ \  & $\left\|u^N - u_h^N\right\|_0$ &$rate_s$ & $\left\|u^{2N} - u_h^{2N}\right\|_0$ &$rate_s$ &  $\left\|u^{3N} - u_h^{3N}\right\|_0$&$rate_s$ \\
\hline
1/8&1.3786e-02&&4.0332e-02&&1.0297e-01&  \\
1/16&3.4693e-03&1.990&1.0142e-02& 1.991&2.5891e-02&1.991  \\
1/32&8.6928e-04&1.996&2.5273e-03&2.004&6.4395e-03&2.007  \\
1/64&2.1796e-04&1.996&6.1931e-04&2.028&1.5652e-03&2.040  \\
\hline
\end{tabular}}
\end{center}
\end{table}
\begin{table}
\caption{Numerical temporal accuracy for different $\alpha.$}\label{t3}
%\resizebox{\textwidth}{40mm}
\begin{center}
\setlength{\tabcolsep}{2.5mm}{
\begin{tabular}{llllllllll}
\hline
$\alpha$& $N$  & $k=1$&   &  $k=2$&  & $k=3$
\\ &  & $E(h,N,k)$ &$rate_t$ & $E(h,N,k)$ &$rate_t$ &  $E(h,N,k)$&$rate_t$ \\
\hline
0.4
&$400$ &4.1103e-03&&3.8838e-04& &8.9101e-04&     \\
&$800$&3.1522e-03&0.382&1.9445e-04&0.998&4.4578e-04&0.999\\
&$1600$ &2.4100e-03&0.387 &9.7313e-05&0.998&2.2299e-04&0.999   \\
&$3200$ &1.8385e-03&0.390&4.8685e-05&0.999&1.1153e-04&0.999    \\
\hline
0.6
&$400$ &1.2744e-03&&4.7951e-04& &1.0567e-03&     \\
&$800$&8.2696e-04&0.623&2.4074e-04&0.994&5.2949e-04&0.996    \\
&$1600$ &5.3841e-04&0.619 &1.2074e-04&0.995&2.6518e-04&0.997   \\
&$3200$ &3.5159e-04&0.614&6.0513e-05&0.996&1.3275e-04&0.998    \\
\hline
\end{tabular}}
\end{center}
\end{table}
\begin{table}
\caption{The temporal accuracy and spatial accuracy for $\alpha=0.4$ at $t_{kN}.$}\label{t4}
%\resizebox{\textwidth}{55mm}
\begin{center}
\setlength{\tabcolsep}{2.5mm}{
\begin{tabular}{llllllllll}
\hline
&$k=1$&& $k=2$ &&  $k=3$& \\
\hline
$N$\ \  & $\left\|u^N - u_h^N\right\|_0$ &$rate_t$ & $\left\|u^{2N} - u_h^{2N}\right\|_0$ &$rate_t$ &  $\left\|u^{3N} - u_h^{3N}\right\|_0$&$rate_t$ \\
\hline
400&1.3858e-04&&3.8838e-04&&8.9101e-04&  \\
800&6.9385e-05&0.998&1.9445e-04&0.990&4.4578e-04&0.999  \\
1600&3.4722e-05&0.998&9.7313e-05&0.998&2.2299e-04&0.999  \\
3200&1.7371e-05& 0.999&4.8685e-05&0.999&1.1153e-04&0.999  \\
\hline
$h$\ \  & $\left\|u^N - u_h^N\right\|_0$ &$rate_s$ & $\left\|u^{2N} - u_h^{2N}\right\|_0$ &$rate_s$ &  $\left\|u^{3N} - u_h^{3N}\right\|_0$&$rate_s$ \\
\hline
1/8&9.7302e-03&&1.9471e-02&&4.5101e-02&  \\
1/16&2.4393e-03&1.996&4.8819e-03&1.995&1.1308e-02 &1.995  \\
1/32&6.1026e-04&1.998&1.2213e-03&1.999&2.8290e-03&1.998  \\
1/64&1.5259e-04&1.999&3.0539e-04&1.999&7.0739e-04&1.999  \\
\hline
\end{tabular}}
\end{center}
\end{table}

\section{Conclusion}\label{Sec6}
In this paper, we derive the exact solution of a time-fractional
diffusion equation with constant time delay, investigate its regularity,
and decompose the solution into its singular and regular parts.
Unlike the previous works which show that the solutions of those
fractional differential equations
with constant time delay have multi-sigularity at $t=k\tau^+$, the solution
we consider is only exhibits singularity
at $t = 0^+$ for its first time derivative, and at both $t = 0^+$ and $\tau^+$
for its second time derivative.
Therefore, we choose the Gr\" unwald-Letnikov type approximation to handle
the time-fractional operators,
and
establish the fully discrete finite element scheme. For the developed numerical scheme,
it is proved to be stable. Then,
in order to discuss the
convergence, the temporal local truncation error based on the decomposition of the solution
mentioned in Theorem 2.1 is discussed. Under a established new discrete Gronwall
inequality, we obtain the error estimate of the proposed numerical scheme. Finally, some numerical tests are provided to
verify our theoretical results.

\section*{Acknowledgements}
The first researcher is supported by the Research Foundation of Education Commission of Hunan Province of China (Grant No. 23A0126) and the 111 Project (Grant No. D23017),
and the fourth researcher is supported by the National Natural Science Foundation of China (Grant No. 12301533).
%the Natural Science Foundation of Hunan Province of China (Nos. ),


\begin{thebibliography}{00}
\addcontentsline{toc}{section}{References}

\bibitem{Erneux2009applied} T. Erneux, Applied Delay Differential Equations, Springer, New York, 2009.

\bibitem{Smith2011An}H. Smith, An Introduction to Delay Differential Equations with Applications to the Life Sciences, Springer, New York, 2011.

\bibitem{Driver1977ordinary} R.D. Driver, Ordinary and Delay Differential Equations, Springer, New York, 1977.

\bibitem{Gopalsamy1992Stability}K. Gopalsamy, Stability and Oscillations in Delay Differential Equations of Population Dynamics, Kluwer Academic Publishers, Dordrecht, 1992 .

\bibitem{Bellen2003numerical}A. Bellen, M. Zennaro, Numerical Methods for Delay Differential Equations, Clarendon Press, New York, 2003.

\bibitem{Balachandran2009delay} B. Balachandran, T. Kalmar-Nagy, D.E. Gilsinn, Delay Differential Equations: Recent Advances and New Directions, Springer, Berlin, 2009.

\bibitem{Carpinteri2014Fractals} A. Carpinteri, F. Mainardi, Fractals and Fractional Calculus in Continuum Mechanics, Springer, New York, 1997.

\bibitem{Hilfer2000Applications} R. Hilfer, Applications of Fractional Calculus in Physics, World Scientific, Singapore, 2000.

\bibitem{Kilbas2006Theory} A.A. Kilbas, H.M. Srivastava, J.J. Trujillo, Theory and Applications of Fractional Differential Equations, Elsevier,  Amsterdam, 2006.

\bibitem{Sun2018A} H. Sun, Y. Zhang, D. Baleanu, W. Chen, Y. Chen, A new collection of real world applications of fractional calculus in science and engineering, Commun. Nonlinear Sci. Numer. Simulat. 64 (2018) 213--231.

\bibitem{Bu2024Finite}W. Bu, S. Guan, X. Xu, Y. Tang, Finite element method for a generalized constant delay diffusion equation, Commun. Nonlinear Sci. Numer. Simulat. 134 (2024) 108015.

\bibitem{Prakash2020Exact}P. Prakash, S. Choudhary, V. Daftardar-Gejji, Exact solutions of generalized nonlinear time-fractional reaction-diffusion equations with time delay, Eur. Phys. J. Plus. 135 (2020) 1--24.

\bibitem{Zhu2016Local}B. Zhu, L. Liu, Y. Wu, Local and global existence of mild solutions for a class of nonlinear fractional reaction-diffusion equations with delay, Appl. Math. Lett. 61 (2016) 73--79.

\bibitem{Li2021Monotone}Q. Li, G. Wang, M. Wei, Monotone iterative technique for time-space fractional diffusion equations involving delay, Nonlinear Anal. Model. Control. 26 (2021) 241--258.

\bibitem{Yao2023Stability}Z. Yao, Z. Yang, Stability and asymptotics for fractional delay diffusion-wave equations, Math. Methods Appl. Sci. 46 (2023) 15208--15225.

\bibitem{Li2018convergence}L. Li, B.  Zhou, X. Chen, Z. Wang, Convergence and stability of compact finite difference method for nonlinear time fractional reaction-diffusion equations with delay, Appl. Math. Comput. 337 (2018) 144--152.

\bibitem{zhao2018fast}Y. Zhao, P. Zhu, W. Luo, A fast second-order implicit scheme for non-linear time-space fractional diffusion equation with time delay and drift term, Appl. Math. Comput. 336 (2018) 231--248.

\bibitem{Li2018ane}T. Li, Q. Zhang, W. Niazi, Y. Xu, M. Ran,  An effective algorithm for delay fractional convection-diffusion wave equation based on reversible exponential recovery method, IEEE Access 7 (2018) 5554--5563.

\bibitem{Zhang2023numerical}Y. Zhang, Z. Wang, Numerical simulation for time-fractional diffusion-wave equations with time delay, J. Appl. Math. Comput. 69 (2023) 137--157.

\bibitem{Alikhanov2024second}A.A. Alikhanov, M.S. Asl, C. Huang, A. Khibiev, A second-order difference scheme for the nonlinear time-fractional diffusion-wave equation with generalized memory kernel in the presence of time delay, J. Comput. Appl. Math. 438 (2024) 115515.

\bibitem{Pimenov2017on}V.G. Pimenov, A.S. Hendy, R.H. De Staelen, On a class of non-linear delay distributed order fractional diffusion equations, J. Comput. Appl. Math. 318 (2017) 433--443.

\bibitem{peng2023convergence}S. Peng, M. Li, Y. Zhao,  F. Wang, Y. Shi, Convergence and superconvergence analysis for nonlinear delay reaction-diffusion system with nonconforming finite element, Numer. Meth. Part. D. E. 39 (2023) 716--743.

\bibitem{Peng2024uncond}S. Peng, M. Li, Y. Zhao, F. Liu, F, Cao, Unconditionally convergent and superconvergent finite element method for nonlinear time-fractional parabolic equations with distributed delay, Numer. Algorithms 95 (2024) 1643--1714.

\bibitem{zaky2023L1}M.A. Zaky,  K. Van Bockstal, T.R Taha, D. Suragand, A.S. Hendy, An L1 type difference/Galerkin spectral scheme for variable-order time-fractional nonlinear diffusion-reaction equations with fixed delay, J. Comput. Appl. Math. 420 (2023) 114832.

\bibitem{Hendy2019novel}A.S. Hendy, J.E. Mac\'ias-D\'iaz, A novel discrete Gr\"onwall inequality in the analysis of difference schemes for time-fractional multi-delayed diffusion equations, Commun. Nonlinear Sci. Numer. Simulat. 73 (2019) 110--119.

\bibitem{Tan2022L1}T. Tan, W. Bu, A. Xiao, L1 method on nonuniform meshes for linear time-fractional diffusion equations with constant time delay, J. Sci. Comput. 92 (2022) 98.

\bibitem{cen2023the}D. Cen, S. Vong, The tracking of derivative discontinuities for delay fractional equations based on fitted L1 method, Comput. Methods Appl. Math.
 23 (2023) 591--601.

\bibitem{cen2023corrected}D. Cen, C.  Ou, S. Vong, Corrected L-type method for multi-singularity problems arising from delay fractional equations, J. Sci. Comput. 97 (2023) 15.

\bibitem{Ou2024Variable}C. Ou, D. Cen, S. Vong, Variable-step L1 method combined with time two-grid algorithm for multi-singularity problems arising from two-dimensional nonlinear delay fractional equations, Commun. Nonlinear Sci. Numer. Simulat. 139 (2024) 108270.

\bibitem{Sakamoto2011Initial}K. Sakamoto, M. Yamamoto, Initial value/boundary value problems for fractional diffusion-wave equations and applications to some inverse problems, J. Math. Anal. Appl. 382 (2011) 426--447.

\bibitem{Luchko2012Initial}Y. Luchko, Initial-boundary-value problems for the one-dimensional time-fractional diffusion equation, Fract. Calc. Appl. Anal. 15 (2012) 141--60.

\bibitem{li2019Theory}C. Li, M. Cai, Theory and Numerical Approximations of Fractional Integrals and Derivatives, SIAM, Philadelphia, 2019.

\bibitem{Lubich1986Discretized}C. Lubich, Discretized fractional calculus, SlAM J. Math. Anal. 17 (1986) 704--719.

\bibitem{Lubich1988Convolution}C. Lubich, Convolution quadrature and discretized operational calculus. I, Numer. Math. 52 (1988) 129--145.

\bibitem{Chen2022Sharp}H. Chen, Y. Shi, J. Zhang, Y. Zhao, Sharp error estimate of a Gr\"unwald-Letnikov scheme for reaction-subdiffusion equations, Numer. Algorithms 89 (2022) 1465--1477.

\bibitem{Liao2019A}H. Liao, W. McLean, J. Zhang, A discrete Gr\"onwall inequality with applications to numerical schemes for subdiffusion problems, SIAM J. Numer. Anal. 57 (2019) 218--237.

\bibitem{Zeng2017Second}F. Zeng, Z. Zhang, G.E. Karniadakis, Second-order numerical methods for multi-term fractional differential equations: smooth and non-smooth solutions, Comput. Methods Appl. Mech. Eng. 327 (2017) 478--502.

\bibitem{Zhao2015A}L. Zhao, W. Deng, A series of high-order quasi-compact schemes for space fractional diffusion equations based on the superconvergent approximations for fractional derivatives, Numer. Meth. Part. D. E. 31 (2015) 1345--1381.

\bibitem{Thomee2006Galerkin}V. Thom\'ee, Galerkin Finite Element Methods for Parabolic Problems, Springer, Berlin, 2006.


\end{thebibliography}
\end{document}